\documentclass{amsart}

\usepackage{amssymb}
\usepackage{hyperref}
\usepackage{graphicx}
\usepackage[text={14.7cm,21.2cm},centering]{geometry}

\newtheorem{defi}{Definition}[section]
\newtheorem{prop}[defi]{Proposition}
\newtheorem{thm}[defi]{Theorem}
\newtheorem{lem}[defi]{Lemma}
\newtheorem{cor}[defi]{Corollary}

\numberwithin{equation}{section}

\newcommand{\N}{\mathbb{N}}
\newcommand{\Z}{\mathbb{Z}}
\newcommand{\R}{\mathbb{R}}
\newcommand{\Haus}{\mathcal{H}}
\newcommand{\D}{\mathcal{D}}
\newcommand{\Mass}{\operatorname{\mathbf M}}
\newcommand{\Norm}{\operatorname{\mathbf N}}
\newcommand{\Int}{\operatorname{\mathbf I}}
\newcommand{\Flat}{\operatorname{\mathbf F}}
\newcommand{\IFlat}{\operatorname{\mathcal F}}
\newcommand{\BV}{\operatorname{BV}}
\newcommand{\Le}{\mathcal L}
\newcommand{\loc}{\text{loc}}
\newcommand{\Lip}{\operatorname{Lip}}
\newcommand{\Lipc}{\Lip_c}
\newcommand{\Lipl}{\Lip_{\loc}}
\newcommand{\Hol}{\operatorname{H}}
\newcommand{\Holc}{\Hol_c}
\newcommand{\Holl}{\Hol_{\loc}}
\newcommand{\spt}{\operatorname{spt}}
\newcommand{\diam}{\operatorname{diam}}
\newcommand{\defl}{\mathrel{\mathop:}=}
\newcommand{\defr}{=\mathrel{\mathop:}}

\begin{document}

\title{Integration of H\"older forms and currents in \linebreak snowflake spaces}
\author{Roger Z\"ust}
\address{Department of Mathematics, ETH Z\"urich, Switzerland}
\email{roger.zuest@math.ethz.ch}
%\date{February 3, 2009}
\thanks{Partially supported by the Swiss National Science Foundation.}
%\subjclass[2000]{49Q15, 26B20}

\begin{abstract}
For an oriented $n$-dimensional Lipschitz manifold $M$ we give meaning to the integral $\int_M f \, dg_1 \wedge \dots \wedge dg_n$ in case the functions $f, g_1, \dots, g_n$ are merely H\"older continuous of a certain order by extending the construction of the Riemann-Stieltjes integral to higher dimensions. More generally, we show that for $\alpha \in (\tfrac{n}{n+1},1]$ the $n$-dimensional locally normal currents in a locally compact metric space $(X,d)$ represent a subspace of the $n$-dimensional currents in $(X, d^\alpha)$. On the other hand, for $n \geq 1$ and $\alpha \leq \tfrac{n}{n+1}$ the vector space of $n$-dimensional currents in $(X, d^\alpha)$ is zero.
\end{abstract}

\maketitle

%\tableofcontents

\section{Introduction}

If $f,g_1, \dots, g_n$ are smooth functions on $[0,1]^n$, the differential form $f \, dg_1 \wedge \dots \wedge dg_n$ is defined and we can calculate the integral
\[ \int_{[0,1]^n} f \, dg_1 \wedge \dots \wedge dg_n. \]
In general, this differential form makes no sense if $f,g_1, \dots, g_n$ are not smooth but merely H\"older continuous. Nevertheless, we want to show that in case the sum of the H\"older exponents of these $n+1$ functions is bigger than $n$, the integral above can be given a reasonable value by generalizing the construction of the classical Riemann-Stieltjes integral to higher dimensions. More precisely, this integral, we call it $\int_{[0,1]^n}f \, d(g_1,\dots,g_n)$, is constructed recursively on the dimension of the cube by approximating it with Riemannian sums of the form
\[ \sum_{i=1}^{2^{kn}}f(\mu_{B_i})\int_{\partial B_i} g_1 \, d(g_2,\dots,g_n), \]
where $B_1, \dots, B_{2^{kn}}$ is the  partition of $[0,1]^n$ into $2^{kn}$ cubes of equal size and $\mu_{B_i}$ is the barycenter of $B_i$. Our result for $n=1$ is covered in \cite{young} where L.C.\ Young showed that the Riemann-Stieltjes integral exists even under weaker assumptions. With the usual partition of unity construction this integral extends to oriented Lipschitz manifolds and a variant of Stokes' theorem for H\"older continuous functions is presented.

In the last section we discuss the connection to the theory of currents in metric spaces. Metric currents have been introduced by Ambrosio and Kirchheim in \cite{kirch}, extending the classical Federer-Fleming theory of \cite{flem} to complete metric spaces. We will mainly work with the local currents introduced by Lang in \cite{lang} not relying on a finite mass assumption. For a locally compact metric space $(X,d)$ the $n$-dimensional currents $\D_n(X)$ are functions
\[ T : \Lipc(X)\times \prod_{i=1}^n \Lipl (X) \rightarrow \R \]
that are $(n+1)$-linear, continuous in a suitable sense and satisfy $T(f,\pi_1,\dots,\pi_n) = 0$ whenever some $\pi_i$ is constant on a neighborhood of the support of $f$. For $\alpha \in (0,1)$ we are interested in the vector space $\mathcal D_n(X,d^\alpha)$ of $n$-dimensional currents in the snowflake space $(X,d^\alpha)$. By approximating H\"older with Lipschitz functions we generalize the result obtained for the Riemann-Stieltjes integral above and show that any locally normal current $T \in \Norm_{n,\loc}(X)$, as defined in Chapter $5$ of \cite{lang}, has a natural extension to a functional 
\[ \bar T: \Holc^\alpha(X)\times \Holl^{\beta_1}(X)\times \dots \times \Holl^{\beta_n}(X) \rightarrow \R \]
on a tuple of H\"older functions if the exponents satisfy $\alpha + \beta_1 + \dots + \beta_n > n$. In particular, if $\alpha = \beta_1 = \dots = \beta_n > \tfrac{n}{n+1}$, this extension is a current in $\mathcal D_n(X,d^\alpha)$ and hence $\Norm_{n,\loc}(X)$ can be identified with a subspace of $\mathcal D_n(X,d^\alpha)$. On the other hand, we show that $\D_n(X,d^\alpha) = \{0\}$ if $n \geq 1$ and $\alpha \leq \tfrac{n}{n+1}$.

\vspace{0.5cm}
{\parindent0mm
\textbf{Acknowledgements:}
I would like to thank Urs Lang for many inspiring discussions and for carefully reading earlier versions of this paper. I am also grateful to Christian Riedweg for some helpful comments.
}

\section{Approximation of H\"older continuous functions}

A map $f$ from $(X,d_X)$ to $(Y,d_Y)$ is said to be H\"older continuous of order $\alpha \in (0,1]$ if there exists a $C \in [0,\infty)$ such that
\[ d_Y(f(x),f(x')) \leq C\,d_X(x,x')^\alpha \]
holds for all $x,x' \in X$. The smallest $C$ with this property is denoted by $\Hol^\alpha(f)$. The set of all such maps is $\Hol^\alpha(X,Y)$ or $\Hol^\alpha(X)$ in case $(Y,d_Y) = (\R,|.|)$. If $\alpha = 1$, we speak of Lipschitz continuous maps and write $\Lip$ instead of $\Hol^1$. If $X$ is bounded, a basic property of these sets is that $\Hol^\beta(X,Y) \subset \Hol^\alpha(X,Y)$ for $0 < \alpha \leq \beta \leq 1$. With respect to the usual addition and multiplication of functions $\Hol^\alpha(X)$ is a vector space, and an algebra if $X$ is bounded. The H\"older exponents are multiplicative with respect to compositions, i.e.\ $g \circ f \in \Hol^{\alpha\beta}(X,Z)$ if $f \in \Hol^\alpha(X,Y)$ and $g \in \Hol^\beta(Y,Z)$.

The next results show how H\"older functions can be approximated by Lipschitz functions. Up to minor modifications of the second lemma they are contained in the appendix of \cite{gromov} written by Stephen Semmes. For $A\subset X$ we denote by $A_\epsilon$ the closed $\epsilon$-neighborhood $\{x\in X \,:\, d(x,A) \leq \epsilon\}$ of $A$.

\begin{lem}\cite[Theorem B.6.3]{gromov}
\label{gromprop}
Let $k > 0$ and $0 < \alpha < 1$ be constants. If $(f_j)_{j\in\Z}$ is a family of functions from $X$ to $\R$ such that
\begin{enumerate}
	\item $\|f_j\|_\infty \leq k2^{j\alpha}$,
	\item $f_j$ is $k2^{j(\alpha-1)}$-Lipschitz,
	\item $\sum_{j\in\Z}f_j(x_0)$ converges for some $x_0\in X$,
\end{enumerate}
then $\sum_{j\in\Z}f_j$ converges uniformly on bounded subsets of $X$ to a function which is H\"older continuous of order $\alpha$ and the partial sums have bounded $\alpha$-H\"older constants. Conversely, every H\"older continuous function of order $\alpha$ admits such a representation.
\end{lem}

\begin{lem}
\label{approx}
Let $C>0$ and $\mathcal F \subset \Hol^\alpha(X)$ be such that $\Hol^\alpha(f) \leq C$ holds for all $f \in \mathcal F$. Then for every $\epsilon > 0$ and $f \in \mathcal F$ we can assign a function $f_\epsilon$ such that
\begin{enumerate}
	\item \label{approx1} $\Lip(f_\epsilon) \leq C\epsilon^{\alpha-1}$,
	\item \label{approx2} $\|f - f_\epsilon\|_\infty \leq C \epsilon^\alpha$,
	\item \label{approx3} $\spt(f_\epsilon) \subset \spt(f)_\epsilon$,
	\item \label{approx4} $\Hol^\alpha(f_\epsilon) \leq 3C$,
	\item \label{approx5} $\|g_\epsilon - h_\epsilon\|_\infty \leq \|g-h\|_\infty$ for all $g,h \in \mathcal F$.
\end{enumerate}
\end{lem}
{\parindent0mm
The following proof is for the most part contained in the proof of \cite[Theorem B.6.16]{gromov}.
}
\begin{proof}
We define $f_\epsilon$ by
\[ f_\epsilon(x) \defl \inf \{f(y) + C\epsilon^{\alpha - 1}d(x,y) \, : \, y \in X\}. \]
$f_\epsilon$ is the infimum of $C\epsilon^{\alpha - 1}$-Lipschitz functions and if finite it is $C\epsilon^{\alpha - 1}$-Lipschitz too. Clearly $f_\epsilon(x) \leq f(x)$. If $d(x,y) \geq \epsilon$, then
\[ f(y) + C\epsilon^{\alpha - 1}d(x,y) \geq f(y) + \Hol^\alpha(f) d(x,y)^\alpha \geq f(x) \]
and therefore 
\[ f_\epsilon(x) = \inf \{f(y) + C\epsilon^{\alpha - 1}d(x,y) \, : \, y \in X, \, d(x,y) \leq \epsilon \}. \]
By this characterization \eqref{approx3} is obvious. In addition
\[ f_\epsilon(x) \leq f(x) \leq \inf \{f(y) + C\epsilon^\alpha : \, y \in X, \, d(x,y) \leq \epsilon \} \leq f_\epsilon(x) + C\epsilon^\alpha \]
for all $x \in X$ which shows \eqref{approx2}. Hence $f_\epsilon(x)$ is finite and as a consequence \eqref{approx1} holds. If $d(x,y)\leq \epsilon$, then
\[ |f_\epsilon(x)-f_\epsilon(y)| \leq C\epsilon^{\alpha-1}d(x,y) \leq Cd(x,y)^\alpha. \]
On the other hand if $d(x,y) \geq \epsilon$ we combine \eqref{approx1} and \eqref{approx2} to conclude \eqref{approx4}:
\begin{align*}
|f_\epsilon(x)-f_\epsilon(y)| & \leq 2\|f-f_\epsilon\|_\infty + |f(x)-f(y)| \\
 & \leq 2C\epsilon^\alpha + C d(x,y)^\alpha \\
 & \leq 2Cd(x,y)^\alpha + Cd(x,y)^\alpha.
\end{align*}
To verify \eqref{approx5} let $g,h \in \mathcal F$. By a straightforward evaluation
\begin{align*}
g_\epsilon(x) & = \inf \{g(y) + C\epsilon^{\alpha - 1}d(x,y) \, : \, y \in X\} \\
 & \leq \|g - h\|_\infty + \inf \{h(y) + C\epsilon^{\alpha - 1}d(x,y) \, : \, y \in X\} \\
 & = \|g - h\|_\infty + h_\epsilon(x)
\end{align*}
and hence \eqref{approx5} holds.
\end{proof}

\section{A generalized Riemann-Stieltjes integral}

\subsection{Construction}

Let $f \in \Hol^\alpha(A)$, $g_1 \in \Hol^{\beta_1}(A), \dots, g_n \in \Hol^{\beta_n}(A)$ be H\"older continuous functions on a box $A = [u_1,v_1] \times \dots \times [u_n,v_n] \subset \R^n$. In this section we define a value for $\int_A f \, d(g_1, \dots, g_n)$, or shorter $\int_A f \, dg$, where $g \defl (g_1, \dots, g_n)$.

The construction of the integral is done recursively. In dimension $0$ the integral is defined to be the evaluation functional. Assuming that the integral in dimensions $0,\dots,n-1$ is already constructed we use the boundary integrals
\[ \int_{\partial B} g_1 \, d(g_2,\dots,g_n) \]
of boxes $B = [s_1,t_1] \times \dots \times [s_n,t_n] \subset A$ to build up the Riemannian sums. They are defined by
\[ \sum_{i=1}^n\sum_{j=0}^1 (-1)^{i+j}\int_{B_{(i,j)}}g_1 \, d(g_2,\dots,g_n), \]
where 
\[ B_{(i,j)} \defl [s_1,t_1] \times \dots \times [s_{i-1},t_{i-1}] \times \{s_i + j(t_i-s_i)\} \times [s_{i+1},t_{i+1}] \times \dots \times [s_n,t_n] \]
and the functions are restricted to these codimension one boxes (to be precise, $B_{(i,j)}$ is identified with a box in $\R^{n-1}$ by omitting the $i$-th coordinate and each function is rearranged accordingly). This is the standard orientation convention for the boundary as used for example in \cite{spivak}. If $B$ is the interval $[s,t]$, the boundary integral is just $\int_{\partial [s,t]} g = g(t) - g(s)$. We define the Riemannian sums
\[ I_n(f,g,\mathcal P,\xi) \defl \sum_{B \in \mathcal P} f(\xi_B) \int_{\partial B} g_1 \, d(g_2,\dots,g_n), \]
where $\mathcal P$ is a partition of $A$ into finitely many boxes with disjoint interiors and $\xi = \{\xi_B\}_{B \in \mathcal P}$ is a collection of points such that $\xi_B \in B$. $I_n(f,g,\mathcal P)$ is the sum above where each $\xi_B$ is assumed to be $\mu_B$, the barycenter of $B$. The mesh, $\|\mathcal P\|$, of a partition $\mathcal P$ of is the maximal diameter of a box in $\mathcal P$. In dimension $1$ this is the usual construction of the Riemann-Stieltjes integral $\int_s^t f \, dg$. It is defined to be the limit, as the mesh of the partition $\mathcal P$ of the interval $[s,t]$ approaches zero, of the approximating sum $I_1(f, g, \mathcal P, \xi)$. To calculate the integral over $A$ we will use only very special partitions, namely $\mathcal P_0(A),\mathcal P_1(A),\mathcal P_2(A),\dots$, where $\mathcal P_0(A)$ consists of the box $A$ alone and $\mathcal P_{k+1}(A)$ is constructed from $\mathcal P_k(A)$ by dividing each box into $2^n$ similar boxes half the size. If it is clear which box is partitioned, we simply write $\mathcal P_k$ instead of $\mathcal P_k(A)$. The definition of $I_n(f,g,\mathcal P,\xi)$ is motivated by Stokes' Theorem. We will make use of it in the following form:

\begin{lem}
\label{boxstokes}
Let $g_1,\dots, g_n$ be Lipschitz functions defined on a box $A \subset \R^n$. Then 
\[ \int_A \det D(g_1,\dots,g_n) \, d\Le^n = \sum_{i=1}^n\sum_{j=0}^1 (-1)^{i+j} \int_{A_{(i,j)}} g_1 \det D(g_2,\dots,g_n) \, d\Le^{n-1}. \]
\end{lem}

For smooth functions the proof is standard and will be omitted. Formulated with differential forms it can be found for example in \cite{spivak}. A Lipschitz function on $\R^n$ can be approximated uniformly by a sequence of smooth functions with bounded Lipschitz constants, see e.g.\ \cite[4.1.2]{fed}, and the integrals in the lemma agree by a continuity argument, see e.g.\ \cite[Example 3.2]{kirch} or \cite[Proposition 2.6]{lang} and the references there for more details.

\begin{thm}
\label{thm2}
For all $n \in \N$, all boxes $A \subset \R^n$ and all numbers $\alpha,\beta_1,\dots,\beta_n$ contained in $(0,1]$ such that $\alpha + \beta_1 + \dots + \beta_n > n$ the function
\[ \int_A: \Hol^\alpha(A) \times \Hol^{\beta_1}(A) \times \dots \times \Hol^{\beta_n}(A) \rightarrow \R \]
\[ (f,g_1,\dots,g_n) \mapsto \int_A f \, dg \defl \lim_{k \rightarrow \infty} I_n(f,g,\mathcal P_k,\xi_k) \]
is well defined and independent of the choice of $(\xi_k)_{k\in\N}$. $\int_A$ satisfies the following properties:
\begin{enumerate}
	\item \label{lin} $\int_A$ is $(n+1)$-linear.
	\item	\label{lip} In case $\beta_1 = \dots = \beta_n = 1$, the identity
	 \[ \int_A f \, dg = \int_A f \det Dg \, d\Le^n \]
	  holds (for smooth $f$ and $g$ this agrees with $\int_A f \, dg_1\wedge \dots \wedge dg_n$).
	\item	\label{cont} $\int_A$ is continuous in the sense that 
\[ \int_A f_m \, dg_m \rightarrow \int_A f \, dg, \; \mbox{ for } m\rightarrow \infty \]
whenever $(f_m)_{m\in\N}$ and $(g_{m,i})_{m\in\N}$ are sequences converging uniformly to $f$ resp.\ $g_i$ on $A$ and $\Hol^\alpha(f_m)$ resp.\ $\Hol^{\beta_i}(g_{m,i})$ are bounded in $m$ for all $i = 1, \dots, n$.
\end{enumerate}
Moreover $\int_A$ is uniquely defined by \eqref{lip} and \eqref{cont}.
\end{thm}

\begin{proof}
Uniqueness of the integral is a direct consequence of Lemma~\ref{approx} (or Lemma~\ref{gromprop}). Every H\"older continuous function can be approximated by Lipschitz functions in such a way that \eqref{cont} applies and by \eqref{lip} the integral for Lipschitz functions is given.

For $n=0$ the theorem is clear. Let $A \subset \R^n$, $f$ and $g$ be as in the theorem and assume that in dimensions $0,\dots,n-1$ the integral is already constructed. It should be noticed that $\beta_1+\dots+\beta_n > n-1$ since $\alpha \leq 1$. So $\int_{\partial B} g_1 \, d(g_2,\dots,g_n)$ does indeed exist for any $n$-dimensional box $B\subset A$ and hence the Riemannian sums $I_n(f,g,\mathcal P_k,\xi_k)$ are well defined. In dimensions $1\leq m \leq n-1$ we additionally assume the existence of constants $C_m'(\hat \beta_1, \dots, \hat\beta_m)$ and $C_m(\hat \alpha, \hat \beta_1, \dots, \hat\beta_m)$ such that for any box $\hat B \subset \R^m$ and functions $\hat f \in \Hol^{\hat \alpha}(\hat B)$, $\hat g_1 \in \Hol^{\hat \beta_1}(\hat B), \dots, \Hol^{\hat \beta_m}(\hat B)$ of orders satisfying $\hat \alpha + \hat \beta_1 + \dots + \hat \beta_m > m$ the following estimates hold:

%depending only on $(H^{\hat \beta_1}, \dots, H^{\hat \beta_m},\hat \beta_1,\dots,\hat \beta_m)$ and $C_m$ depending only on $(H^{\hat \alpha},H^{\hat \beta_1}, \dots, H^{\hat \beta_m},\hat \alpha,\hat \beta_1, \dots, \hat \beta_m)$ such that 
\begin{align}
\label{inducprop2}
\left|\int_{\partial \hat B} \hat g_1 \, d(\hat g_2,\dots,\hat g_m)\right| & \leq C_m' \diam(\hat B)^{\sum_{i=1}^m \hat \beta_i} \prod_{i = 1}^m \Hol^{\hat \beta_i}(\hat g_i), \\
\label{inducprop}
\left|\int_{\hat B}\hat f \, d\hat g - I_m(\hat f,\hat g,\mathcal P_k(\hat B),\hat \xi_k)\right| & \leq C_m \diam(\hat B)^{\hat \alpha + \sum_{i=1}^m \hat \beta_i} \\
\nonumber
 & \quad \; \cdot 2^{k(m - \hat \alpha - \sum_{i=1}^m \hat \beta_i)} \Hol^{\hat \alpha}(\hat f) \prod_{i = 1}^m \Hol^{\hat \beta_i}(\hat g_i).
\end{align}

We will show the existence of $C_n'$ and $C_n$ and these two estimates for the box $A$. Let $B \subset A$ be any $n$-dimensional box. Firstly, an estimate for the boundary integral
\[ J(B,g) \defl \int_{\partial B} g_1 \, d(g_2,\dots,g_n) \]
is established. For $n=1$, \eqref{inducprop2} holds with $C_1'(\beta) \defl 1$ because
\[ \left|\int_{\partial [s,t]} g\right| = |g(t) - g(s)| \leq \Hol^\beta(g)(t-s)^\beta. \]
Now let $n > 1$. To shorten notation we define $\bar \beta \defl \sum_{i=1}^n \beta_i$, $\gamma \defl \alpha + \bar \beta$, $H^\beta \defl \prod_{i = 1}^n \Hol^{\beta_i}(g_i)$ and $H^{\alpha,\beta} \defl \Hol^\alpha(f) H^\beta$. Setting $k=0$ and $m=n-1$ in \eqref{inducprop} leads to
\[ \left| J(B,g) - \sum_{i=1}^n\sum_{j=0}^1(-1)^{i+j} g_1(\xi_{B_{(i,j)}}) \int_{\partial B_{(i,j)}} g_2 \, d(g_3,\dots,g_n) \right| \leq 2n C_{n-1}(\beta) \diam(B)^{\bar \beta} H^\beta. \]
The identity
\[ \sum_{i=1}^n\sum_{j=0}^1(-1)^{i+j} \int_{\partial B_{(i,j)}}g_2 \, d(g_3, \dots, g_n) = 0 \]
is true because of the orientation convention and applying \eqref{inducprop2} with $m=n-1$ to the faces $B_{(i,j)}$ results in
\begin{align}
\nonumber
\left| J(B,g) \right| & \leq \sum_{i=1}^n\sum_{j=0}^1 |g_1(\xi_{B_{(i,j)}}) - g_1(\xi_B)| \left|\int_{\partial B_{(i,j)}}g_2 \, d(g_3,\dots,g_n) \right| \\
\nonumber
 & \quad \; + \, 2n C_{n-1}(\beta) \diam(B)^{\bar \beta} H^\beta \\
\nonumber
 & \leq 2n \Hol^{\beta_1}(g_1) \diam(B)^{\beta_1} C_{n-1}'(\beta_2,\dots,\beta_n) \diam(B)^{\sum_{i = 2}^m \beta_i} \prod_{i = 2}^m \Hol^{\beta_i}(g_i) \\
\nonumber
 & \quad \; + \, 2n C_{n-1}(\beta) \diam(B)^{\bar \beta} H^\beta \\
\label{eq2}
 & = C'_n(\beta) \diam(B)^{\bar \beta} H^\beta,
\end{align}
where
\begin{equation*}
%\label{ind3}
C'_n(\beta) \defl 2n \left(C_{n-1}'(\beta_2,\dots,\beta_n) + C_{n-1}(\beta) \right).
\end{equation*}
Next we show that the limit $\int_A f \, dg$ exists. To do this we will make use of the identity
\begin{equation}
\label{cancel}
J(B,g) = \sum_{i=1}^n\sum_{j=0}^1 (-1)^{i+j} \sum_{F \in \mathcal P_1(B_{(i,j)})} \int_{F} g_1 \, d(g_2,\dots,g_n) = \sum_{\tilde B \in \mathcal P_1(B)} J(\tilde B,g).
\end{equation}
This is true because the integrals over faces of some $\tilde B$ which are not contained in $\partial B$ cancel in pairs. Now we determine an upper bound for $|I_n(f,g,\mathcal P_{k+1}) - I_n(f,g,\mathcal P_k)|$. If $B \in \mathcal P_k$ for some $k \geq 1$, denote by $B'$ the box in $\mathcal P_{k-1}$ with $B \subset B'$. By \eqref{cancel}
\begin{align}
\nonumber
\left|I_n(f,g,\mathcal P_k) - I_n(f,g,\mathcal P_{k-1})\right| & = \left|\sum_{B \in \mathcal P_k} f(\mu_B) J(B,g) - f(\mu_{B'}) \sum_{B \in \mathcal P_k} J(B,g)\right| \\
\nonumber
 & \leq \sum_{B \in \mathcal P_k} |f(\mu_B)-f(\mu_{B'})| \left|J(B,g) \right| \\
\nonumber
 & \leq \Hol^\alpha(f) \frac{\diam(A)^\alpha}{2^{k\alpha}} \sum_{B \in \mathcal P_k} \left|J(B,g) \right| \\
\nonumber
 & \leq 2^{kn} \Hol^\alpha(f) \frac{\diam(A)^\alpha}{2^{k\alpha}} C'_n(\beta) \frac{\diam(A)^{\bar \beta}}{2^{k \bar \beta}} H^\beta \\
\label{eq3}
 & = C'_n(\beta)\diam(A)^\gamma2^{k(n-\gamma)} H^{\alpha,\beta}.
\end{align}
The last inequality holds by \eqref{eq2} and the fact that the cardinality of $\mathcal P_k$ is exactly $2^{kn}$. Hence $(I_n(f,g,\mathcal P_k))_{k\in\N}$ is a Cauchy sequence because $\gamma > n$ and the limit
\[ \int_A f \, dg = \lim_{k\rightarrow \infty}I_n(f,g,\mathcal P_k) \]
exists as stated in the theorem. To additionally handle the intermediate points we note that
\begin{align*}
\left|I_n(f,g,\mathcal P_k,\xi_k) - I_n(f,g,\mathcal P_k)\right| & \leq \sum_{B \in \mathcal P_k} |f(\xi_B) - f(\mu_B)| \left|J(B,g)\right| \\
 & \leq C'_n(\beta) \diam(A)^\gamma 2^{k(n-\gamma)} H^{\alpha, \beta}
\end{align*}
analogously to the estimate above. To show \eqref{inducprop} in case $m=n$ we calculate
\begin{align}
\nonumber
\left|\int_A f \, dg - I_n(f,g,\mathcal P_k,\xi_k)\right| \leq & \, \left|I_n(f,g,\mathcal P_k) - I_n(f,g,\mathcal P_k,\xi_k)\right| \\
\nonumber
& \; + \sum_{j=k+1}^\infty\left|I_n(f,g,\mathcal P_j) - I_n(f,g,\mathcal P_{j-1})\right| \\
\nonumber
 \leq & \, C'_n(\beta)\diam(A)^\gamma H^{\alpha,\beta} \sum_{j=k}^\infty 2^{j(n-\gamma)} \\
\label{eq2b}
 = & \, C_n(\alpha,\beta) \diam(A)^\gamma 2^{k(n-\gamma)} H^{\alpha,\beta},
\end{align}
where
\[ C_n(\alpha,\beta) \defl \frac{C'_n(\beta)}{1 - 2^{n-\gamma}}. \]

We now prove the remaining parts of the theorem. In dimension $0$ the integral is linear and by induction it is multilinear in all dimensions because the approximating sums are. The proof of \eqref{lip} is by induction on the dimension $n$. Let $f \in \Hol^\alpha(A)$ have $\alpha$-H\"older constant $H$ and the $g_i \in \Lip(A)$, $i=1,\dots,n$, have $L$ as a common Lipschitz constant. If $\mathcal P$ is a partition of $A$, then
\begin{align*}
\left|I_n(f,g,\mathcal P) - \int_A f(x) \det Dg(x) \, d\Le^n(x) \right| & = \left|\sum_{B \in \mathcal P}\int_B (f(\mu_B) - f(x)) \det Dg(x) \, d\Le^n(x) \right| \\
 & \leq \sum_{B \in \mathcal P}\int_B |f(\mu_B) - f(x)| |\det Dg(x)| \, d\Le^n(x) \\
 & \leq HL^n \sum_{B \in \mathcal P}\int_B \|\mu_B - x\|^\alpha \, d\Le^n(x) \\
 & \leq HL^n \|\mathcal P\|^\alpha \Le^n(A)
\end{align*}
which is small if $\|\mathcal P\|$ is small. The first equation needs justification. If $n=1$, it holds because
\[ \int_s^t g' \, d\mathcal{L} = g(t) - g(s) = \int_{\partial [s,t]} g \]
by basic analysis or by Lemma~\ref{boxstokes} applied to an interval $[s,t] \subset A$. If $n > 1$,
\begin{align*}
\int_{B} \det Dg \, d \Le^{n} & = \sum_{i=1}^n\sum_{j=0}^1 (-1)^{i+j} \int_{B_{(i,j)}} g_1 \det D(g_2, \dots, g_n) \, d \Le^{n-1} \\
 & = \sum_{i=1}^n\sum_{j=0}^1 (-1)^{i+j} \int_{B_{(i,j)}} g_1 \, d(g_2, \dots, g_n) \\
 & = \int_{\partial B} g_1 \, d(g_2, \dots, g_n)
\end{align*}
and these equations are true on any $n$-dimensional box $B \subset A$ successively by Lemma~\ref{boxstokes} on $B$, the induction hypothesis and the definition of $\int_{\partial B}$.

The remaining part is the proof of \eqref{cont}, the continuity of $\int_A$. Let $f_m$ and $g_m$ be two sequences with the described properties and let $H$ be a common upper bound for the H\"older constants of all functions involved. By the multiliearity of $\int_A$ it suffices to consider bounded sequences $f_m,g_{m,1},\dots,g_{m,n}$ one of which converges to zero and to conclude that the integral $\int_A f_m \, dg_m$ converges to zero as well. Combining \eqref{eq2b} and \eqref{eq2}
\begin{align}
\nonumber
\left|\int_{A}f_m \, dg_m\right| & \leq |I_n(f_m,g_m,\mathcal P_k,\xi_k)| +
 C_n(\alpha,\beta) \diam(A)^{\gamma}2^{k(n - \gamma)}H^{n+1} \\
\nonumber
 & \leq \|f_m\|_\infty \sum_{B \in \mathcal P_k}\left|\int_{\partial B} g_m\right| + C_n(\alpha,\beta) \diam(A)^{\gamma}2^{k(n - \gamma)} H^{n+1} \\
\label{dddd}
 & \leq C'_n(\beta) \|f_m\|_\infty \diam(A)^{\bar \beta}2^{k(n - \bar \beta)} H^n + C_n(\alpha,\beta)\diam(A)^{\gamma}2^{k(n - \gamma)} H^{n+1}.
\end{align}
Using either the third estimate if $f_m$ tends to zero or the second and the induction hypothesis if $g_{m,i}$ tends to zero for some $i$ leads to
\[ \limsup_{m\rightarrow \infty} \left|\int_{A}f_m \, dg_m\right| \leq C_n(\alpha,\beta)\diam(A)^{\gamma}2^{k(n - \gamma)} H^{n+1}. \]
By varying $k$ this expression is arbitrary small. This completes the proof of the theorem.
\end{proof}

For $n=1$ the theorem restates some results of \cite{young} by L.C.\ Young. It is shown there that in the setting of the theorem, the Riemann-Stieltjes integral $\int_s^t f \, dg$ exists (thereby allowing all partitions, not only the special ones we used) and has all the properties described. More generally, the Riemann-Stieltjes integral over $[s,t]$ is defined for functions $f \in W_p([s,t])$ and $g \in W_q([s,t])$ if they have no common discontinuities and $1/p + 1/q > 1$. $W_p([s,t])$ denotes the space of functions $f: [s,t] \rightarrow \R$ with bounded mean variation of order $p$, i.e.
\[ \sup \sum_{i=1}^m |f(x_i) - f(x_{i-1})|^p < \infty, \]
where the supremum ranges over all partitions $s = x_0 \leq x_1 \leq \dots \leq x_m = t$ and all $m \in \N$. If $p\geq 1$, the class $W_p([s,t])$ contains the H\"older functions $\Hol^{1/p}([s,t])$.

It is rather easy to compute $\int_A f \, dg$ numerically. If we know only the values of $f$ and $g$ on the corners of the boxes in $\mathcal P_k(A)$ for some $k$, we can recursively compute approximations for the integral on lower-dimensional sub-boxes to get an approximation of $\int_A f \, dg$ in the end. In case $g$ is Lipschitz, this enables us to numerically calculate $\int_A f \det Dg \, d\Le^n$ without taking $\det Dg$ into account, which, by the way, might not be defined on a set of measure zero.

\subsection{Necessity of the assumption on the H\"older exponents}
\label{assump}

The next example demonstrates that the bound on the H\"older exponents is sharp if we want a function $\int_A$ as in the theorem satisfying \eqref{lip} and \eqref{cont}. Let $\alpha,\beta_1,\dots,\beta_n$ be real numbers in the interval $(0,1]$ such that $\gamma = \alpha + \beta_1 + \dots + \beta_n \leq n$. We consider the box $A = [0,2\pi]^n$ and smooth functions $f_m, g_{m,1}, \dots, g_{m,n}$ defined by
\begin{align*}
f_m(x) & \defl \sum_{i=1}^m \frac{1}{2^{i\alpha}} \sin(2^i x_1) \dots \sin(2^i x_n), \\
g_{m,k}(x) & \defl \sum_{i_k=1}^m \frac{1}{2^{i_k\beta_k}} \cos(2^{i_k}x_k), \; \mbox{ for } k = 1,\dots,n.
\end{align*}
By Lemma~\ref{gromprop} these functions converge to H\"older continuous functions $f$, $g_1, \dots, g_n$ in a way \eqref{cont} of Theorem~\ref{thm2} is applicable and by \eqref{lip} we can calculate $\int_A f_m\, dg_m$. Because $D g_m$ is diagonal as a matrix with respect to the standard basis of $\R^n$, the integrand is given by
\[ f_m \det Dg_m = \sum_{i,i_1,\dots,i_n} \left( 2^{-i\alpha + \sum_{k=1}^n i_k(1-\beta_k)} \prod_{k=1}^n \sin(2^i x_k)\sin(2^{i_k} x_k) \right). \]
For $l,l'\in \N$ it holds that $\int_0^{2\pi} \sin(lx) \sin(l'x) \, dx$ is $\pi$ if $l=l'$ and $0$ otherwise. This identity together with Fubini's theorem implies that only the summands with $i = i_1 = \dots = i_n$ need to be considered. Hence
\begin{align*}
\int_A f_m \det Dg_m \, d\Le^n & = \sum_{i=1}^m 2^{-i\alpha + \sum_{k=1}^n i(1-\beta_k)} \prod_{k=1}^n \int_0^{2\pi} \sin^2(2^ix_k) \, dx_k \\
 & = \pi^n \sum_{i=1}^m 2^{i(n - \gamma)} .
\end{align*}
This sum is unbounded in $m$. Consequently, $\int_A$ can't be extended to include the functions $f,g_1,\dots,g_n$.

\subsection{Some Properties}

We now list some properties of the integral defined above. The proofs will be rather short and rely on approximation by Lipschitz functions, Lemma~\ref{approx}, and the properties of the integral in Theorem~\ref{thm2}.

\begin{prop}
\label{properties}
Let $A$,$f$ and $g$ be given as in Theorem~\ref{thm2}.
\begin{enumerate}
	\item \label{prop1} (additivity) If $\mathcal P$ is a partition of $A$, then $\int_A f \, dg = \sum_{A_i\in \mathcal P} \int_{A_i} f \, dg$.
	\item \label{prop2} $f$ is extended to be zero outside of $A$. If this extension is continuous on a box $B \supset A$ and each $g_i$ is extended arbitrarily to a H\"older continuous function of the same order, then $\int_A f \, dg = \int_B f \, dg$.
	\item \label{prop3} (locality) If some function $g_i$ is constant on a neighborhood of the support of $f$, then $\int_A f \, dg = 0$.
	\item \label{prop4} (alternating property) If $g_i = g_j$ for some different $i$ and $j$, then $\int_A f \, dg = 0$.
	\item \label{prop5} (product rule) If $\alpha \leq \beta_1$ and $h,h'\in \Hol^{\beta_1}(A)$, then
	\[ \int_A f \, d(hh',g_2,\dots,g_n) = \int_A fh \, d(h',g_2,\dots,g_n) + \int_A fh' \,  d(h,g_2,\dots,g_n). \]
	\item \label{prop6} (parametrization property) Let $U$ be a neighborhood of $\spt(f)$ in $A$, $\varphi: U \rightarrow \varphi(U) \subset \R^n$ a bi-Lipschitz map with $\det D \varphi \geq 0$ $\Le^n$-almost everywhere and $B\supset \varphi(U)$ a box such that $\varphi(U \cap \partial A) = \varphi(U) \cap \partial B$. If $f \circ \varphi^{-1}$ is extended to be zero on $B \setminus \varphi(U)$ and each function $g_i \circ \varphi^{-1}$ is extended arbitrarily to a H\"older continuous function of the same order, then
	\[ \int_A f \, d(g_1, \dots, g_n) = \int_B f \circ \varphi^{-1} \, d(g_1 \circ \varphi^{-1}, \dots, g_n \circ \varphi^{-1}). \]
\end{enumerate}
\end{prop}

\begin{proof}
By approximation \eqref{prop1},\eqref{prop3} and \eqref{prop4} are direct consequences of the respective results for Lipschitz functions. The same is true for \eqref{prop5} observing that for $\psi \in \Hol^\alpha(A)$, $\psi' \in \Hol^\beta(A)$ and $\alpha \leq \beta$, the product $\psi \psi'$ satisfies $\Hol^\alpha(\psi \psi') \leq \|\psi\|_\infty\Hol^\beta(\psi')\diam(A)^{\beta - \alpha} + \|\psi'\|_\infty\Hol^\alpha(\psi)$.

If the extension of $f$ in \eqref{prop2} is continuous, it is clearly H\"older continuous of the same order as $f$. Therefore \eqref{prop2} is a consequence of \eqref{prop1}.

The integral on the right-hand side of \eqref{prop6} is well defined because the extension of $f \circ \varphi^{-1}$ is H\"older continuous of the same order as $f$ and the same holds for each $g_i \circ \varphi^{-1}$. For Lipschitz functions the identity is a special case of the area formula, see e.g.\ \cite[3.2.3]{fed}. The general result follows by approximating $f$, $g_1, \dots, g_n$ by Lipschitz functions as before, with the addition that every approximation of $f$ has support in $U$.
\end{proof}

Next we give an upper bound for $|\int_A f \, dg|$ similar to \eqref{dddd} which takes the shape of $A$ a bit better into account.

\begin{cor}
\label{estimcor}
Let $A$,$f$ and $g$ be as in Theorem~\ref{thm2} and let $\epsilon > 0$ be the length of the shortest edge of $A$. Then 
\begin{align}
\nonumber
\left|\int_A f \, dg \right| & \leq K_n \left(\|f\|_\infty \epsilon^{\bar \beta - n} + \Hol^\alpha(f) \epsilon^{\gamma - n} \right) \mathcal{L}^n(A) \prod_{i=1}^n \Hol^{\beta_i}(g_i) \\
\nonumber
 & \leq K_n' \left(\|f\|_\infty \epsilon^{\bar \beta - (n - 1)} \Haus^{n-1}(\partial A) + \Hol^\alpha(f) \epsilon^{\gamma - n} \mathcal{L}^n(A) \right) \prod_{i=1}^n \Hol^{\beta_i}(g_i)
\end{align}
for $\bar \beta = \sum_{i=1}^n \beta_i$, $\gamma = \alpha + \bar \beta$ as before and some constants $K_n$ and $K_n'$ depending only on $n$,$\alpha$ and $\beta = (\beta_1, \dots, \beta_n)$.
\end{cor}

\begin{proof}
Let $n \geq 2$. The edge of $A$ parallel to the $i$-th coordinate axis has length $x_i$ and w.l.o.g.\ $\epsilon = x_1$. For $j=2,\dots,n$ we write $x_j = (N_j + \delta_j)\epsilon$ where $N_j \in \N$ and $\delta_j \in [0,1)$. We partition $A$ into $N_2 \cdot N_3 \cdots N_n$ cubes with edge length $\epsilon$ and some smaller boxes. Applying \eqref{eq2} and \eqref{eq2b} with $k=0$  to these cubes gives in combination with \eqref{prop1} of Proposition~\ref{properties}
\begin{align*}
\left|\int_A f \, dg \right| & \leq \left( C_n'(\beta) \|f\|_\infty (\sqrt n \epsilon)^{\bar \beta} + C_n(\alpha,\beta) \Hol^\alpha(f) (\sqrt n \epsilon)^{\gamma} \right) \prod_{i=1}^n \Hol^{\beta_i}(g_i) \prod_{j=2}^n (N_j + 1) \\
 & \leq K_n'' \left( \|f\|_\infty \epsilon^{\bar \beta} + \Hol^\alpha(f) \epsilon^\gamma \right) \prod_{i=1}^n \Hol^{\beta_i}(g_i) \prod_{j=2}^n 2\frac{x_j}{\epsilon}
\end{align*}
where $K_n'' = \tfrac{K_n}{2^{n-1}}$ and the first inequality of the corollary is immediate. The second is a direct consequence of the first by noting that $2n \Le^n(A) \leq \epsilon \Haus^{n-1}(\partial A)$. The case $n=1$ is clear since $\epsilon = \diam(A) = \Le(A)$.
\end{proof}

\subsection{Stokes' theorem for H\"older continuous functions}

The integral of Theorem~\ref{thm2} satisfies 
\begin{equation}
\label{stokesimp}
\int_A 1 \, d(g_1, \dots, g_n) = \int_{\partial A} g_1 \, d(g_2, \dots, g_n)
\end{equation}
by definition. The goal here is to extend this integral to oriented Lipschitz manifolds and to show that a similar variant of Stokes' theorem holds in this setting.

A metric space $(M,d)$ is said to be an $n$-dimensional Lipschitz manifold if it can be covered by charts $(U,\varphi)$, where $U$ is an open bounded subset of $M$ and $\varphi$ is a bi-Lipschitz map of $U$ onto an open bounded subset of $\{x \in \R^n : x_1 \leq 0\}$. The boundary $\partial M$ of $M$ is the set of those points that are mapped into $\{x \in \R^n : x_1 = 0\}$ by some (and hence all) charts. $\partial M$ is either empty or an $(n-1)$-dimensional Lipschitz manifold (with the induced metric). If $X$ is a paracompact Hausdorff space with an atlas of charts whose transition functions are bi-Lipschitz, a metric $d$ can be constructed such that $(X,d)$ is a Lipschitz manifold compatible with this atlas. This and related results can be found in \cite{luukk}.

In the comments below we assume that $n \geq 2$. The $1$-dimensional manifolds (with their $0$-dimensional boundaries) need special considerations and are left to the reader. $M$ is said to be orientable if there exists an atlas $\{(U_i,\varphi_i)\}_{i\in I}$ such that every transition function $\varphi_i\circ {\varphi_k}^{-1}$ is orientation preserving in the sense that $\det D (\varphi_i\circ {\varphi_k}^{-1})$ is positive $\mathcal{L}^n$-almost everywhere on $\varphi_k(U_i\cap U_k)$. 
%Because there is a constant $L \geq 1$ such that both $\varphi_i \circ {\varphi_k}^{-1}$ and $\varphi_k\circ {\varphi_i}^{-1}$ are $L$-Lipschitz, this implies that $\det D(\varphi_i\circ {\varphi_k}^{-1}) \geq 1/L^n$ $\mathcal{L}^n$-almost everywhere on $\varphi_k(U_i\cap U_k)$. An orientation is a maximal atlas of such charts.
An orientation on $M$ induces an orientation on $\partial M$ with the defining property that for every positively oriented chart $\varphi : U \rightarrow \{x_1 \leq 0\}$ of $M$, the restriction $\varphi|_{\partial M}$ is a positively oriented chart of $\partial M$ (given the obvious identification of $\R^{n-1}$ with $\{ x \in \R^n : x_1 = 0\}$ by deleting the first coordinate).
%Because $\{x_1 = 0\} \subset \R^n$ is a set of measure zero it is not totally obvious that such an orientation exists. It can be established by applying Lemma~\ref{boxstokes} to a transition function $\psi$ of two positively oriented charts with domain $V \subset \{x_1 \leq 0\}$ and any box $A_\epsilon$ with one face on $V \cap \{x_1 = 0\}$, called $F$, and the opposite side on $\{x_1 = -\epsilon\}$. Here $\epsilon > 0$ is assumed small enough such that $A_\epsilon$ is contained in $V$. By letting $\epsilon$ tend to zero it is possible to conclude that
%\begin{equation*}
%\int_F \det D(\psi|_{\{x_1 = 0\}}) \, d\mathcal{L}^{n-1} \geq 0.
%\end{equation*}
%This is true for any such box $F$ and therefore $\det D(\psi|_{\{x_1 = 0\}})$ is non-negative $\mathcal{L}^{n-1}$-almost everywhere on $V \cap \{x_1 = 0\}$ by the Lebesgue differentiation theorem.

The next result extends the integral of the last section to oriented Lipschitz manifolds and states a variant of Stokes' theorem.

\begin{thm}
\label{stokethm}
\begin{sloppypar}
Let $M$ be an oriented $n$-dimensional Lipschitz manifold and let $\alpha$, $\beta_1,\dots,\beta_n$ be constants contained in the interval $(0,1]$ such that $\alpha + \beta_1 + \dots + \beta_n > n.$ Then there is a unique multilinear function
\end{sloppypar}
\[ \int_M: \Holc^\alpha(M)\times \Holl^{\beta_1}(M)\times \dots \times \Holl^{\beta_n}(M) \rightarrow \R, \quad (f,g_1,\dots,g_n) \mapsto \int_M f \, d(g_1,\dots,dg_n) \]
such that 
\[ \int_M f \, d(g_1,\dots,g_n) = \int_B f \circ \varphi^{-1} \, d(g_1\circ \varphi^{-1},\dots,g_n\circ \varphi^{-1}) \]
whenever $(U,\varphi)$ is a positvely oriented chart which contains $\spt(f)$ and $B$ is a box with $\varphi(U) \subset B \subset \{x_1 \leq 0\}$. Furthermore, if $f = 1$ on a neighborhood of $\spt(g_1)$ (which has to be compact for this reason), then
\begin{equation}
\label{aaa2}
\int_M f \, d(g_1,\dots,g_n) = \int_{\partial M} g_1 \, d(g_2,\dots,g_n).
\end{equation}
\end{thm}

\begin{proof}
Let $\{(U_i,{\varphi_i})\}_{1\leq i \leq N}$ be finitely many positively oriented charts such that the $U_i$ cover $\mbox{spt}(f)$. We choose a Lipschitz partition of unity $\{\theta_i\}_{1\leq i \leq N}$ for $\mbox{spt}(f)$ subordinate to this covering. Assuming the multilinearity and parametrization property the integral has to be defined by
\[ \int_M f \, d(g_1,\dots,g_n) \defl \sum_{i=1}^N \int_{{\varphi_i}(U_i)} (\theta_i f)\circ {\varphi_i}^{-1} \, d\left(g_1\circ {\varphi_i}^{-1}, \dots, g_n\circ {\varphi_i}^{-1}\right). \]
$\int_{{\varphi_i}(U_i)}$ means an integral as defined in Theorem~\ref{thm2} over a box in $\{x_1 \leq 0\}$ that contains ${\varphi_i}(U_i)$ and each $g_j \circ \varphi_i^{-1}$ is extended arbitrarily to a H\"older continuous function of the same order. The right-hand side is well defined by \eqref{prop2} and \eqref{prop3} of Proposition~\ref{properties}. It is independent of the charts and the partition of unity by \eqref{prop6} of Proposition~\ref{properties}. $\int_M$ is multilinear because the $\int_{{\varphi_i}(U_i)}$ are.

To show Stokes' theorem for this integral we consider this time a cover $\{(U_i,{\varphi_i})\}_{1\leq i \leq N}$ of $\spt(g_1)$ and a subordinate Lipschitz partition of unity $\{\theta_i\}_{1\leq i \leq N}$ for $\spt(g_1)$. By \eqref{prop3} of Proposition~\ref{properties} we can replace $f$ by a function which is $1$ on every $U_i$ without changing the left-hand side of \eqref{aaa2}. W.l.o.g.\ $U_1,\dots, U_{N'}$ are those sets that meet the boundary $\partial M$. By the linearity of $\int_M$ in the second argument and \eqref{stokesimp}
\begin{align*}
\int_M f \, d(g_1,\dots,g_n) & = \sum_{i=1}^N \int_M f \, d(\theta_i g_1, g_2,\dots, g_n) \\
 & = \sum_{i=1}^N \int_{{\varphi_i}(U_i)} d\left((\theta_i g_1)\circ {\varphi_i}^{-1}, g_2\circ {\varphi_i}^{-1}, \dots, g_n\circ {\varphi_i}^{-1}\right) \\
 & = \sum_{i=1}^{N'} \int_{{\varphi_i}(U_i) \cap \{x_1 = 0\}} (\theta_i g_1)\circ {\varphi_i}^{-1} \,  d\left(g_2\circ {\varphi_i}^{-1}, \dots, g_n\circ {\varphi_i}^{-1}\right) \\
 & = \int_{\partial M} g_1 \, d(g_2,\dots,g_n).
\end{align*}
\end{proof}

\section{Currents in snowflake spaces and applications}

\subsection{Preliminaries}

Following \cite{lang} the vector space $\D_n(X)$ of $n$-dimensional currents in a locally compact metric space $(X,d)$ are those functions $T: \D^n(X) \rightarrow \R$, where 
\[ \D^n(X) = \Lipc(X)\times \prod_{i=1}^n \Lipl (X), \]
such that:
\begin{enumerate}
	\item (multilinearity) $T$ is $(n+1)$-linear.
	\item (locality) If $n \geq 1$, $T(f, \pi^1, \dots, \pi^n) = 0$ whenever some $\pi^i$ is constant on a neighborhood of $\spt(f)$.
	\item (continuity) $T$ is continuous in the sense that
\[ T(f_j, \pi_j) \rightarrow T(f, \pi), \; \mbox{ for } j\rightarrow \infty, \]
whenever $(f_j)_{j\in\N}$ and $(\pi^i_{j})_{j\in\N}$, $i = 1,\dots n$, are sequences which satisfy the following convergence criteria:
\begin{enumerate}
	\item $f_j$ converges uniformly to $f$, the Lipschitz constants $\Lip(f_j)$ are bounded in $j$ and there is a compact set which contains $\spt(f_j)$ for all $j$.
	\item For every compact set $K \subset X$ the Lipschitz constants $\Lip(\pi^i_j|_K)$ are bounded in $j$ and $\pi^i_j|_K$ converges uniformly to $\pi^i|_K$.
\end{enumerate}
\end{enumerate}
\vspace{0.2cm}
Here are some definitions related to a current $T \in \D_n(X)$ we will need:

\begin{itemize}
	\item If $n\geq 1$, the \textbf{boundary} of $T$ is the current $\partial T \in \D_{n-1}(X)$ given by
\[ \partial T(f, \pi^1, \dots, \pi^{n-1}) \defl T(\sigma, f, \pi^1, \dots, \pi^{n-1}), \]
where $\sigma \in \Lipc(X)$ is any function such that $\sigma = 1$ on a neighborhood of $\spt(f)$, see \cite[Definition 3.4]{lang}.

	\item The \textbf{support} of $T$ is the smallest closed set $\spt(T) \subset X$ such that $T(f,\pi) = 0$ whenever $\spt(f) \cap \spt(T) = \emptyset$, see \cite[Definition 3.1]{lang}.

	\item Let $Y$ be a locally compact space, $A$ a locally compact subset of $X$ containing $\spt(T)$ and $F \in \Lipl(A,Y)$ a proper map, i.e.\ $F^{-1}(K)$ is compact if $K\subset Y$ is compact. Then the \textbf{pushforward} of $T$ via $F$ is the current $F_{\#}T \in \D_n(Y)$ defined by
\[ F_{\#}T(f, \pi^1, \dots, \pi^n) \defl T_A(f \circ F, \pi^1 \circ F, \dots, \pi^n \circ F) \]
for $(f, \pi^1, \dots, \pi^n) \in \D^n(Y)$, see \cite[Definition 3.6]{lang}. $T_A$ denotes the restriction of $T$ to $\D^n(A)$.

	\item The \textbf{mass} of $T$ on an open set $V\subset X$, $\Mass_V(T)$, is the least number $M \in [0,\infty]$ such that 
\[ \sum_{\lambda \in \Lambda} T(f_\lambda,\pi_\lambda) \leq M \]
whenever $\Lambda$ is a finite set, $\spt(f_\lambda)\subset V$, $\sum_{\lambda \in \Lambda}|f_\lambda| \leq 1$ and $\pi^i_\lambda$ is $1$-Lipschitz for all $i$ and $\lambda$, see \cite[Definition 4.2]{lang}. We set $\Norm_V(T) \defl \Mass_V(T) + \Mass_V(\partial T)$ if $n \geq 1$ and $\Norm_V(T) \defl \Mass_V(T)$ if $n=0$. If $V = X$, the index in $\Mass_V$ and $\Norm_V$ is omitted.
\end{itemize}
In $\D_n(X)$ the following subspaces are of special interest: \\
currents with finite mass 
\[ \Mass_n(X) \defl \{T \in \D_n(X) \, : \, \Mass(T) < \infty \} \]
currents with locally finite mass
\[ \Mass_{n,\loc}(X) \defl \{T \in \D_n(X) \, : \, \Mass_V(T) < \infty \text{ for all open } V \Subset X \} \]
normal currents
\[ \Norm_n(X) \defl \{T \in \D_n(X) \, : \, \Norm(T) < \infty \} \]
locally normal currents
\[ \Norm_{n,\loc}(X) \defl \{T \in \D_n(X) \, : \, \Norm_V(T) < \infty \text{ for all open } V \Subset X \} \]

For any $\alpha \in (0,1)$ a snowflake space $(X,d^\alpha)$ is obtained. By abuse of notation we write $X^\alpha$ if it is clear which metric is meant. Obviously $\Hol^\alpha(X) = \Lip(X^\alpha)$, $\Holc^\alpha(X) = \Lipc(X^\alpha)$ and $\Holl^\alpha(X) = \Lipl(X^\alpha)$. The next result points out some basic facts about currents in snowflake spaces.

\begin{lem}
Let $X$ be a locally compact metric space. Every current in $\D_n(X^\alpha)$ is the unique extension of a current in $\D_n(X)$ and $\Mass_{n,\loc}(X^\alpha) = \{0\}$ for $\alpha \in (0,1)$ and $n\geq 1$.
\end{lem}

\begin{proof}
If $g \in \Lip(B)$, where $B$ is a bounded metric space, then
\begin{equation}
\label{estimhol}
\Hol^\alpha(g) \leq \Lip(g) \diam(B)^{1 - \alpha}.
\end{equation}
From this estimate we infer that $\Lipc(X) \subset \Holc^\alpha(X)$ and $\Lipl(X) \subset \Holl^\alpha(X)$ and a current $T \in \D_n(X^\alpha)$ can be restricted to $\D^n(X)$. This restriction defines a current in $\D_n(X)$. The multilinearity and locality axioms are immediate and the continuity axiom holds by \eqref{estimhol}. On the other hand, $T$ is defined by its values on $\Holc^\alpha(X)^{n+1}$ by the locality property, and these functions can be approximated by Lipschitz function as given in Lemma~\ref{approx}. The continuity property then implies that the restriction of $T$ to $\Lipc(X)^{n+1}$ is enough to reconstruct its values on $\D^n(X^\alpha)$.

In the definition of the mass above it would be equally valid to assume that each $\pi^i_\lambda$ is only locally $1$-Lipschitz. Let $n \geq 1$ and $T \in \D_n(X^\alpha)$ with $T \neq 0$. By definition there is a $(f,\pi) \in \Lipc(X)^{n+1}$ with $\|f\|_\infty \leq 1$ and $T(f,\pi) > 0$. It follows from \eqref{estimhol} that $r\pi^1$ has locally an arbitrary small $\alpha$-H\"older constant for every $r \in \R$ and as a result
\[ \Mass_V(T) \geq T(f,r\pi^1,\pi^2,\dots,\pi^n) = rT(f,\pi) \rightarrow \infty \]
for $r \rightarrow \infty$, where $V$ is any open neighborhood of $\spt(f)$ with compact closure.
\end{proof}

\subsection{Extension of locally normal currents}

In what follows we discuss the question whether a current in $\D_n(X)$ can be extended to a current in $\D_n(X^\alpha)$. From Theorem~\ref{stokethm} it follows that every oriented $n$-dimensional Lipschitz manifold $M$ defines a current in $\D_n(M^\alpha)$ for $\alpha > \tfrac{n}{n+1}$. Further extensions are not possible in general as indicated by the counterexample in subsection~\ref{assump}. There are however currents not extendable this far. Let $a=(a_m)_{m\in\N}$ be a sequence of positive numbers such that $\sum_{m=1}^\infty a_m < \infty$. Let $s_m \defl \sum_{k=1}^m a_k$ be the partial sums and set $s_0 \defl 0$. For such a sequence we denote by $[a]$ the current in $\Mass_1(\R)$ induced by $\bigcup_{m=0}^\infty \left[2s_m, 2s_m + a_{m+1} \right]$. We claim that $[a]$ is a current in $\D_1(\R^\alpha)$ if and only if $\sum_{m=1}^\infty a_m^\alpha$ is finite and $\alpha > \tfrac{1}{2}$. If $[a]$ is a current in $\D_1(\R^\alpha)$, then $\alpha > \tfrac{1}{2}$ by the counterexample in subsection~\ref{assump}. We choose a function $g \in \Holc^\alpha(\R)$ such that $g(x) = (x-2s_m)^\alpha$ on each interval $[2s_m, 2s_m + a_{m+1}]$ for all $m \geq 0$. Because $\partial[a] \in \D_0(\R^\alpha)$ we have
\[ \partial [a](g) = \sum_{m = 1}^\infty a_m^\alpha < \infty. \]
On the other hand if $\alpha > \tfrac{1}{2}$ and this sum is finite, then by Corollary~\ref{estimcor}
\begin{equation}
\label{estimex}
\sum_{m=1}^\infty \int_{2s_{m-1}}^{2s_{m-1} + a_m}f \, dg \leq C \sum_{m=1}^\infty \left(\Hol^\alpha(f) \Hol^\alpha(g) a_m^{2\alpha} + \|f\|_\infty \Hol^\alpha(g)a_m^\alpha \right) < \infty,
\end{equation}
for a constant $C$ depending only on $\alpha$. Consequently $[a] \in \D_1(\R^\alpha)$. %Let $\operatorname{I}_a \subset (0,1]$ be the interval of all $\alpha$ such that $[a] \in \D_1(\R^\alpha)$. For example if $a_m = m^{-\beta^{-1}}$, then $\operatorname{I}_a = \left(\beta,1\right] \cap (\tfrac{1}{2},1]$, or if $a_m = m^{-(1 + (\log_2\log_2 m)^{-1})\beta^{-1}}$ for $m \geq 4$, then $\operatorname{I}_a = \left[\beta,1\right] \cap (\tfrac{1}{2},1]$. Similar examples in higher dimensions exist too. These $[a] \in \D_1(\R)$ are flat chains with finite mass but infinite boundary mass. They are in particular not locally normal and therefore it may still be possible that $\Norm_{n,\loc}(X) \subset \D_n(X^\alpha)$ for $\alpha > \tfrac{n}{n+1}$. This turns out to be true and is implied by the next theorem.
For example if the sequence is $a_m = m^{-\beta^{-1}}$, then $[a] \in \D_1(\R^\alpha)$ exactly if $\alpha \in \left(\beta,1\right] \cap (\tfrac{1}{2},1]$, or if $a_m = m^{-(1 + (\log_2\log_2 m)^{-1})\beta^{-1}}$ for $m \geq 4$, then $[a] \in \D_1(\R^\alpha)$ exactly if $\alpha \in \left[\beta,1\right] \cap (\tfrac{1}{2},1]$. Similar examples in higher dimensions exist too. These $[a] \in \D_1(\R)$ are flat chains with finite mass but infinite boundary mass. They are in particular not locally normal and therefore it may still be possible that $\Norm_{n,\loc}(X) \subset \D_n(X^\alpha)$ for $\alpha > \tfrac{n}{n+1}$. This turns out to be true and is implied by the next theorem.

\begin{thm}
\label{normal}
Let $X$ be a locally compact metric space and let $\alpha,\beta_1,\dots,\beta_n$ be constants contained in the interval $(0,1]$ such that $\alpha + \beta_1 + \dots + \beta_n > n$. Then for every $T\in \Norm_{n,\loc}(X)$, there is a unique extension
\[ \bar T: \Holc^\alpha(X)\times \Holl^{\beta_1}(X)\times \dots \times \Holl^{\beta_n}(X) \rightarrow \R \]
such that:
\begin{enumerate}
	\item \label{normal1} $\bar T$ is $(n+1)$-linear,
	\item \label{normal2} $\bar T(f, \pi^1, \dots, \pi^n) = 0$ if some $\pi^i$ is constant on a neighborhood of $\spt(f)$,
	\item \label{normal3} $\bar T$ is continuous in the sense that
\[ \bar T (f_j,\pi_j) \rightarrow \bar T (f,\pi), \; \mbox{ for } j\rightarrow \infty, \]
whenever $(f_j)_{j\in\N}$ and $(\pi^i_{j})_{j\in\N}$, $i = 1,\dots n$, are sequences which satisfy the following convergence criteria:
\begin{enumerate}
	\item $f_j$ converges uniformly to $f$, the H\"older constants $\Hol^\alpha(f_j)$ are bounded in $j$ and there is a compact set which contains $\spt(f_j)$ for all $j$.
	\item For every compact set $K \subset X$ the H\"older constants $\Hol^{\beta_i}(\pi^i_{j}|_K)$ are bounded in $j$ and $\pi^i_{j}|_K$ converges uniformly to $\pi^i|_K$.
\end{enumerate}
\end{enumerate}
\end{thm}

\begin{proof}
Uniqueness is a consequence of Lemma~\ref{approx}. By \eqref{normal1} and \eqref{normal2} of the theorem we can assume that $\pi^1,\dots,\pi^n$ have support contained in a compact neighborhood of $\spt(f)$. With \eqref{approx3} and \eqref{approx4} of Lemma~\ref{approx} this tuple of functions with compact support can be approximated by Lipschitz functions in such a way that the continuity property of $\bar T$ applies. But the value of $\bar T$ for Lipschitz-tuples with compact support is given. So, if there is such an extension, it is unique.

We first consider the case where
\[ (f,\pi^1,\dots,\pi^n) \in \Holc^\alpha(X)\times \Holc^{\beta_1}(X)\times \dots \times \Holc^{\beta_n}(X) \]
and all these functions have support contained in the compact set $K \subset X$. Let $V$ be any open set containing $K$ with finite $\Norm_V(T)$. The latter condition certainly holds if $V$ has compact closure. Let $C^\alpha \geq \Hol^\alpha(f), C^{\beta_1} \geq \Hol^{\beta_1}(\pi^1), \dots, C^{\beta_n} \geq \Hol^{\beta_n}(\pi^n)$ be constants. If $\delta \geq \epsilon > 0$ are small enough such that $K_\delta$ is contained in $V$, we choose approximations $f_\epsilon, \pi^1_\epsilon, \dots, \pi^n_\epsilon$ and $f_\delta', \pi'^1_\delta, \dots, \pi'^n_\delta$ satisfying \eqref{approx1} and \eqref{approx2} of Lemma~\ref{approx} with respect to the constants above in place of $C$ such that all the approximating functions have compact support contained in $V$. We are interested in a bound on the difference
\[ |T(f'_\delta, \pi'_\delta) - T(f_\epsilon, \pi_\epsilon)|. \]
This term is dominated by the sum
\[ |T(f'_\delta - f_\epsilon, \pi'_\delta)| + \sum_{i=1}^n \left |T(f_\epsilon, \pi^1_\epsilon, \dots, \pi^{i-1}_\epsilon, \pi'^i_\delta - \pi^i_\epsilon, \pi'^{i+1}_\delta, \dots, \pi'^n_\delta) \right|. \]
 To shorten notation we write $\bar \beta \defl \sum_{i=1}^n \beta_i$, $\gamma \defl \alpha + \bar \beta$, $C^\beta \defl \prod_{i=1}^n C^{\beta_i}$ and $C^{\alpha,\beta} \defl C^\alpha C^\beta$. Using \cite[Theorem 4.3(4)]{lang}
\begin{align}
\nonumber
|T(f'_\delta - f_\epsilon, \pi_\delta')| & \leq \Mass_V(T) \|f_\delta' - f_\epsilon\|_\infty \prod_{i=1}^n \Lip(\pi'^i_\delta) \\
\nonumber
 & \leq \Mass_V(T) (\|f'_\delta - f\|_\infty + \|f - f_\epsilon\|_\infty) \prod_{i=1}^n \Lip(\pi'^i_\delta) \\
\nonumber
 & \leq \Mass_V(T) (\delta^\alpha + \epsilon^\alpha) C^\alpha \prod_{i=1}^n C^{\beta_i} \delta^{\beta_i - 1} \\
\label{firstest}
 & \leq \Mass_V(T) (\tfrac{\delta^\alpha}{\epsilon^\alpha} + 1) C^{\alpha,\beta} \epsilon^{\gamma - n}.
\end{align}
Assuming that $\epsilon$ is small enough such that $C^\alpha \epsilon^\alpha \leq \|f\|_\infty$ equation (5.1) in \cite{lang} gives
\begin{align*}
\Mass_V(\partial(T \lfloor f_\epsilon)) & \leq \|f_\epsilon\|_\infty \Mass_V(\partial T) + \Lip(f_\epsilon)\Mass_V(T) \\
 & \leq (\|f\|_\infty + C^\alpha \epsilon^\alpha ) \Mass_V(\partial T) + C^\alpha \epsilon^{\alpha - 1} \Mass_V(T) \\
  & \leq 2 \|f\|_\infty \Mass_V(\partial T) + C^\alpha \epsilon^{\alpha - 1} \Mass_V(T).
\end{align*}
An estimate for the terms 
\begin{align*}
S_i \defl & \left |T(f_\epsilon, \pi^1_\epsilon, \dots, \pi^{i-1}_\epsilon, \pi'^i_\delta - \pi^i_\epsilon, \pi'^{i+1}_\delta, \dots, \pi'^n_\delta) \right| \\
 = & \left| \partial(T \lfloor f_\epsilon)(\pi'^i_\delta - \pi^i_\epsilon, \pi^1_\epsilon, \dots, \pi^{i-1}_\epsilon, \pi'^{i+1}_\delta, \dots, \pi'^n_\delta) \right|
\end{align*}
is given by
\begin{align}
\nonumber
S_i & \leq \Mass_V(\partial(T \lfloor f_\epsilon)) \|\pi'^i_\delta - \pi^i_\epsilon\|_\infty \prod_{j=1}^{i-1} \Lip(\pi^j_\epsilon) \prod_{j=i+1}^{n} \Lip(\pi'^j_\delta) \\
\nonumber
 & \leq \Mass_V (\partial(T \lfloor f_\epsilon)) (\tfrac{\delta^{\beta_i}}{\epsilon^{\beta_i}} + 1) C^\beta \epsilon^{\bar\beta-(n-1)} \prod_{j=i+1}^{n} \tfrac{\delta^{\beta_j-1}}{\epsilon^{\beta_j-1}} \\
\label{secondest}
 & \leq \left( 2 \|f\|_\infty \Mass_V(\partial T) C^\beta \epsilon^{\bar\beta-(n-1)} + \Mass_V(T) C^{\alpha,\beta} \epsilon^{\gamma - n} \right) (\tfrac{\delta^{\beta_i}}{\epsilon^{\beta_i}} + 1).
\end{align}
If in addition $\delta \leq 2 \epsilon$, then combining \eqref{firstest} and \eqref{secondest} leads to
\begin{equation}
\label{thirdeq}
|T(f'_\delta, \pi'_\delta) - T(f_\epsilon, \pi_\epsilon)| \leq 6n \Mass_V(\partial T) \|f\|_\infty C^\beta \epsilon^{\bar\beta - (n-1)} + 3(n+1) \Mass_V(T) C^{\alpha,\beta} \epsilon^{\gamma - n}.
\end{equation}
By assumption $\bar\beta - (n-1) \geq \gamma - n > 0$ and $\Mass_V(T) + \Mass_V(\partial T) = \Norm_V(T) < \infty$ and the estimate above implies that $(T(f_{2^{-m}},\pi_{2^{-m}}))_{m\in\N}$ is a Cauchy sequence in $\R$. $\bar T(f,\pi)$ is defined to be its limit and we show now that it does not depend on the choice of the approximating sequence. If $m$ is big enough such that $K_{2^{-m}}\subset V$ and $C^\alpha 2^{-m\alpha} \leq \|f\|_\infty$, we have
\begin{align}
\nonumber
|\bar T(f,\pi) - T(f_{2^{-m}},\pi_{2^{-m}})| & \leq \sum_{j = m+1}^\infty \left[ 6n \Mass_V(\partial T) \|f\|_\infty C^\beta 2^{j((n-1) - \bar\beta)} \right.\\
\label{estim3}
& \left. \qquad\qquad\;\; + \, 3(n+1) \Mass_V(T) C^{\alpha,\beta} 2^{j(n - \gamma)} \right].
\end{align}
Let $2^{-m} \leq \delta \leq 2^{-(m-1)}$ be such that $K_{\delta}\subset V$. Combining \eqref{thirdeq} and \eqref{estim3} gives
\begin{align}
\nonumber
|\bar T(f,\pi) - T(f'_{\delta},\pi'_{\delta})| & \leq |T(f_{2^{-m}},\pi_{2^{-m}}) - T(f'_{\delta},\pi'_{\delta})| \\
\nonumber
& \quad \; + \sum_{j = m+1}^\infty \left[ 6n \Mass_V(\partial T) \|f\|_\infty C^\beta 2^{j((n-1) - \bar\beta)} \right.\\
\nonumber
& \left. \qquad\qquad\quad\;\;\, + \, 3(n+1) \Mass_V(T) C^{\alpha,\beta} 2^{j(n - \gamma)} \right] \\
\nonumber
& \leq \frac{6n+3}{1 - 2^{n - \gamma}} \left[ \Mass_V(\partial T) \|f\|_\infty  C^\beta 2^{m((n-1) - \bar\beta)} \right. \\
\label{estim4}
& \qquad\qquad\quad\;\;\, \left. \; + \, \Mass_V(T) C^{\alpha,\beta} 2^{m(n - \gamma)} \right].
\end{align}
This shows that $T(f_\epsilon,\pi_\epsilon) \rightarrow \bar T(f,\pi)$ for $\epsilon \rightarrow 0$ whenever the approximating functions $f_\epsilon$, $\pi^1_\epsilon, \dots, \pi^n_\epsilon$ satisfy \eqref{approx1} and \eqref{approx2} of Lemma~\ref{approx} with an upper bound on the constants used in place of $C$ and the supports of these functions are compact and contained in a fixed open set $V$ with $\Norm_V(T) < \infty$.

Next we show that $\bar T$ is linear in the first argument. The other cases are done similarly. Let $g$ be another function in $\Holc^\alpha(X)$ and assume that the support of $g$ is also contained in $K$. To handle the sum $f + g$ we set the approximation $(f+g)_\epsilon'$ to be $f_\epsilon + g_\epsilon$. This is an approximation for $f+g$ such that \eqref{approx1} and \eqref{approx2} of Lemma~\ref{approx} holds with $C = \Hol^\alpha(f) + \Hol^\alpha(g)$ and the support of $(f+g)'_\epsilon$ is contained in $K_\epsilon$. Since $T((f + g)'_\epsilon,\pi_\epsilon) = T(f_\epsilon,\pi_\epsilon) + T(g_\epsilon,\pi_\epsilon)$ we get
\begin{align*}
|\bar T(f + g,\pi) - \bar T(f,\pi) - \bar T(g,\pi)| & \leq |\bar T(f + g,\pi) - T((f+g)_\epsilon',\pi_\epsilon)| \\
& \quad \; + \, |\bar T(f,\pi) - T(f_\epsilon,\pi_\epsilon)| + |\bar T(g,\pi) - T(g_\epsilon,\pi_\epsilon)|,
\end{align*}
where the latter sums tend to zero by \eqref{estim4} if $\epsilon$ tends to zero. Multiplication by a constant is done likewise.

Assume now that $\pi^i$ is constant on a neighborhood of $\spt(f)$, w.l.o.g.\ $i=1$. Let $c$ be the value of $\pi^1$ on $\spt(f)$. The approximation $\pi^1_\epsilon$ as constructed in the proof of Lemma~\ref{approx} satisfies $\spt(\pi^1_\epsilon) \subset \spt(\pi^1)_\epsilon$ and similarly $\spt(\pi^1_\epsilon - c) \subset \spt(\pi^1 - c)_\epsilon$. If $\epsilon$ is small enough such that $\spt(f)_{\epsilon} \cap \spt(\pi^1-c)_{\epsilon} = \emptyset$, then $T(f_\epsilon, \pi_\epsilon) = 0$ and consequently $\bar T(f, \pi) = 0$.

If $(f,\pi)$ is an element of $\Holc^\alpha(X)\times \Holl^{\beta_1}(X)\times \dots \times \Holl^{\beta_n}(X)$, we choose $\varphi \in \Lipc(X)$ such that $\varphi = 1$ on a neighborhood of $\spt(f)$. An easy calculation shows that 
\begin{equation}
\label{holest}
\Hol^{\beta_i}(\varphi \pi^i) \leq \|\varphi\|_\infty\Hol^{\beta_i}(\pi^i|_{\spt(\varphi)}) + \|\pi^i|_{\spt(\varphi)}\|_\infty \Hol^{\beta_i}(\varphi)
\end{equation}
and we can define
\begin{equation}
\label{holest2}
\bar T(f,\pi) \defl \bar T(f, \varphi\pi^1, \dots, \varphi\pi^n).
\end{equation}
By the locality and multilinearity property just proven, this definition does not depend on $\varphi$. It is clear that these two properties also hold on $\Holc^\alpha(X)\times \Holl^{\beta_1}(X)\times \dots \times \Holl^{\beta_n}(X)$.

If $f_j, \pi^1_j,\dots, \pi^n_j$ are sequences as given in the theorem, there is a compact set $K \subset X$ such that $\spt{f_j} \subset K$, $\Hol^\alpha(f_j) \leq H$ and $\|f_j\|_\infty \leq B$ for all $j$. With \eqref{holest} and the definition in \eqref{holest2} we can assume (by maybe enlarging $K$ and $H$) that $\spt(\pi^i_j) \subset K$ and $\Hol^{\beta_i}(\pi^i_j) \leq H$ for all $i$ and $j$. To apply Lemma~\ref{approx} let $\mathcal F$ be the collection of all these functions and set $C=H$. By the multilinearity of $\bar T$, in order to show the convergence of $|\bar T(f,\pi) - \bar T(f_j,\pi_j)|$ to zero, we can assume that one of the sequences $f_j, \pi^1_j, \dots, \pi^n_j$ converges uniformly to zero. If $m$ is big enough such that $H2^{-m\alpha} \leq B$ and $K_{2^{-m}}\subset V$ for some open set $V$ with $\Mass_V(T) < \infty$, we can apply \eqref{estim4} to conclude that
\begin{align*}
\left|\bar T(f_j,\pi_j) \right| & \leq \left| T((f_j)_{2^{-m}},(\pi_j)_{2^{-m}}) \right| \\
& \quad \; + \, \frac{6n+3}{1 - 2^{n - \gamma}} \left[\Mass_V(\partial T) B H^n 2^{m((n-1) - \bar\beta)} + \Mass_V(T) H^{n+1}2^{m(n - \gamma)} \right].
\end{align*}
$T((f_j)_{2^{-m}},(\pi_j)_{2^{-m}}) \rightarrow 0$ for $j \rightarrow \infty$ by \eqref{approx1}, \eqref{approx3} and \eqref{approx5} of Lemma~\ref{approx} and the continuity of $T$. So, there is an $N \in \N$ such that
\[ \limsup_{j\rightarrow \infty} \left|\bar T(f_j,\pi_j) \right| \leq \frac{6n+3}{1 - 2^{n - \gamma}} \left[\Mass_V(\partial T) B H^n 2^{m((n-1) - \bar\beta)} + \Mass_V(T) H^{n+1}2^{m(n - \gamma)} \right] \]
for all $m \geq N$. Therefore $\lim_{j\rightarrow \infty} |\bar T(f_j,\pi_j) | = 0$ and this concludes the proof of the theorem.
\end{proof}

Let $U$ be an open subset of $\R^n$. By \cite[Theorem 7.2]{lang} the locally normal currents $\Norm_{n,\loc}(U)$ can be identified with the space of functions of locally bounded variation $\BV_{\loc}(U)$. This is the space of all $u \in \text{L}^1_{\loc}(U)$ with
\[ \sup \left\{\int_V u \operatorname{div}(\psi) \, d \Le^n \, : \, \psi \in C^1_c(V,\R^n), \, \|\psi\|_\infty \leq 1 \right\} < \infty \]
for all open sets $V \Subset U$. This identification assigns to $u \in \text{BV}_{\loc}(U)$ the current $[u]$ given by
\[ [u](f,\pi) = \int_U u f \det(D\pi)\, d\Le^n \]
for all $(f,\pi) \in \D^n(U)$. The theorem above gives a meaning to this integral in case the functions $f,\pi_1,\dots,\pi_n$ are only H\"older continuous and thereby extends the scope of Theorem~\ref{thm2}, where $u$ is the characteristic function of a box. But compared to the construction in the proof above the generalized Riemann-Stieltjes integral seems to have some advantages. Firstly, it is rather direct to compute numerically and secondly, for thin boxes the upper bounds calculated in Corollary~\ref{estimcor} are stronger and allow for example the estimate \eqref{estimex}.

\subsection{Applications}

Theorem~\ref{normal} extends several constructions that are known to work for Lipschitz maps to some classes of H\"older maps. For example if $A$ is a locally compact subset of $X$, $Y$ is another locally compact metric space and $\varphi$ is a proper map contained in $\Holl^\alpha(A,Y)$ for some $\alpha > \tfrac{n}{n+1}$, then for every $T \in \Norm_{n,\loc}(X)$ with $\spt(T) \subset A$ the pushforward $\varphi_{\#} T \in \D_n(Y)$ exists. Or if $(u, v_1, \dots, v_k) \in \Holl^\alpha(X) \times \Holl^{\beta_1}(X) \times \dots \times \Holl^{\beta_k}(X)$, where $n \geq k \geq 0$, the current $T \lfloor (u, v_1, \dots, v_k) \in \D_{n-k}(X)$ exists if $\alpha + \beta_1 + \dots + \beta_k > k$, see \cite[Definition 2.3]{lang} for the definition of this construction.

To illustrate this, let $[\Omega] = [\R^2] \lfloor \Omega \in \D_2(\R^2)$ be the current representing the von Koch snowflake domain $\Omega$. We want to find a closed expression for $\partial [\Omega]$ not relying on $[\Omega]$. In the usual way $\Omega$ is constructed as the union of closed sets $\Omega_1 \subset \Omega_2 \subset \dots \subset \Omega$. Each $\Omega_i$ is bi-Lipschitz equivalent to the closed unit ball in $\R^2$. So, there are bi-Lipschitz functions $\varphi_i : S^1 \rightarrow \R^2$ such that ${\varphi_i}_{\#}[S^1] = \partial [\Omega_i]$ by Stokes' theorem. The $\varphi_i$'s can be chosen in such a way that they converge uniformly to a function $\varphi$ and $\Hol^\alpha(\varphi_i)$ is bounded in $i$ for $\alpha = \tfrac{\log 3}{\log 4}$, the reciprocal of the Hausdorff dimension of $\partial \Omega$, see e.g.\ \cite[p.151]{tukia}. This in particular implies that $\varphi$ is H\"older continuous of order $\alpha$ and ${\varphi_i}_{\#}[S^1]$ converges weakly to $\varphi_{\#}[S^1]$ as currents in $\D_1(\R^2)$ due to the fact that $\alpha > \tfrac{1}{2}$. Because $\Mass([\Omega] - [\Omega_i]) \rightarrow 0$ the boundaries $\partial[\Omega_i]$ converge weakly to $\partial [\Omega]$ and hence $\partial [\Omega] = \varphi_{\#}[S^1]$. This leads to the expression
\[ \partial [\Omega](f,g) = \int_{S^1} f \circ \varphi \, d(g \circ \varphi) \]
for all $(f,g) \in \D^1(\R^2)$. In this case the pushforward $\varphi_{\#}[S^1]$ is an integral flat chain. For the definition of flat chains, $\Flat_n(V)$, and integral flat chains, $\IFlat_n(V)$, in an open set $V \subset \R^m$ we refer to \cite[4.1.12]{fed} and \cite[4.1.24]{fed}. The next proposition generalizes this observation about the von Koch curve.

\begin{prop}
Let $T \in \Norm_n(X)$ with compact support and $\varphi \in \Hol^\alpha(X,\R^m)$ for some $\alpha > \tfrac{n}{n+1}$. Then $\varphi_\# T \in \D_n(\R^m)$ is a flat chain respectively an integral flat chain if $T \in \Int_n(X)$.
\end{prop}

\begin{proof}
Because $T$ has compact support we can assume that $X = \spt(T)$ is compact. Let $0 \leq a < b \leq 1$. As in Theorem 3.2 of \cite{weng} the functional $T_a^b$ on $\D^{n+1}([0,1] \times X)$ defined by
\[ T_a^b(f,\pi^1,\dots,\pi^{n+1}) \defl \sum_{i=1}^{n+1} (-1)^{i+1} \int_a^b T(f_t \partial_t \pi^i_t, \pi^1_t, \dots, \pi^{i-1}_t,\pi^{i+1}_t, \dots, \pi^{n+1}_t) \, dt \]
is a current in $\Norm_{n+1}([0,1] \times X)$ resp.\ $\Int_{n+1}([0,1] \times X)$ if $T \in \Int_n(X)$. It satisfies
\begin{equation}
\label{bound}
\partial(T_a^b) = T_b - T_a - (\partial T)_a^b,
\end{equation}
where for any $s \in [0,1]$ the current $T_s$ in $\Norm_n([0,1] \times X)$ resp.\ $\Int_n([0,1] \times X)$ is defined by 
\[ T_s(f,\pi^1,\dots,\pi^n) \defl T(f_s,\pi^1_s,\dots,\pi^n_s). \]
Clearly
\begin{equation}
\label{mass}
\Mass(T_a^b) \leq (b-a)\Mass(T).
\end{equation}
Motivated by Lemma~\ref{approx} we define $\tilde \varphi : [0,1] \times X \rightarrow \R^m$ coordinate-wise by
\[ \tilde \varphi^k(t,x) \defl \inf \{\varphi^k(y) + C t^{\alpha - 1}d(x,y) \, : \, y \in X, \, d(x,y) \leq t \}, \quad k=1,\dots,m, \]
where $C \defl \max_{1 \leq k \leq m} \Hol^\alpha(\varphi^k)$. A simple calculation shows that each $\tilde \varphi^k(.,x)$ is $Cts^{\alpha-2}$-Lipschitz on $[s,t]$ for all $0 < s < t \leq 1$. Together with \eqref{approx1} of Lemma~\ref{approx} the map $\tilde \varphi$ is Lipschitz on $[s,1] \times X$ for all $s \in (0,1]$ and there is a constant $C' > 0$ such that $\tilde \varphi$ is $C's^{\alpha-1}$-Lipschitz on $[s,2s] \times X$ for all $s \in (0,\tfrac{1}{2}]$. Hence, with \eqref{mass}
\[ \Mass(\tilde \varphi_\# T_s^{2s}) \leq (C's^{\alpha-1})^{n+1} s \Mass(T) = C'^{n+1} s^{\alpha (n+1) - n} \Mass(T). \]
This shows that $(\tilde \varphi_\# T_{2^{-i}}^1)_{i \in \N}$ is a Cauchy sequence in $\Norm_{n+1}(\R^m)$ resp.\ $\Int_{n+1}(\R^m)$ equipped with the $\Mass$-norm. This sequence converges to a current $\tilde \varphi_\# T_0^1 \in \Mass_{n+1}(\R^m)$ because $\Mass_{n+1}(\R^m)$ equipped with the $\Mass$-norm is a Banach space by \cite[Proposition 4.2]{lang}. By \cite[Theorem 5.5]{lang} the metric mass and the Euclidean mass are comparable, thus $\tilde \varphi_\# T_0^1$ is in $\Flat_{n+1}(\R^m) \cap \Mass_{n+1}(\R^m)$ resp.\ $\mathcal R_{n+1}(\R^m)$ by \cite[4.1.17]{fed} and \cite[4.1.24]{fed}. With \eqref{bound} this shows that
\begin{align*}
\varphi_\# T & = \lim_{i \rightarrow \infty} \tilde \varphi_\#(T_{2^{-i}}) \\
 & = \lim_{i \rightarrow \infty} \tilde \varphi_\# \left(T_1 - (\partial T)_{2^{-i}}^1 - \partial(T_{2^{-i}}^1) \right) \\
 & = \tilde \varphi_\# T_1 - \lim_{i \rightarrow \infty} \tilde \varphi_\# (\partial T)_{2^{-i}}^1 - \lim_{i \rightarrow \infty} \partial (\tilde \varphi_\# T_{2^{-i}}^1) \\
 & = \tilde \varphi_\# T_1 - \tilde \varphi_\# (\partial T)_0^1 - \partial (\tilde \varphi_\# T_0^1),
\end{align*}
and this is a current in $\Flat_n(\R^m)$ resp.\ $\IFlat_n(\R^m)$.
\end{proof}

The following result shows that many locally normal currents can be realized as pushforwards of Euclidean currents. It is a direct consequence of the Assouad embedding theorem.

\begin{cor}
Let $T$ be in $\Norm_{n,\loc}(X)$ resp.\ $\Int_{n,\loc}(X)$ for a locally compact doubling metric space $X$. Then, for any $\alpha \in (\tfrac{n}{n+1},1)$ there is an $m \in \N$, an open subset $U$ of $\R^m$, a current $S$ in $\Flat_n^{\loc}(U)$ resp.\ $\IFlat_n^{\loc}(U)$ and a bi-Lipschitz map $\varphi: \spt(S) \rightarrow \spt(T)^\alpha$ such that $T = \varphi_{\#} S$. Note that in particular $\varphi \in \Lipl(\spt(S),X)$.
\end{cor}

\begin{proof}
By the Assouad embedding theorem (see e.g.\ \cite{assouad}) there is a bi-Lipschitz embedding $\psi : \spt(T)^\alpha \rightarrow \R^m$ for some $m \in \N$. The image $\psi(\spt(T))$ is a locally compact subset of $\R^m$. By a characterization of such sets $\psi(\spt(T)) = U \cap A$, where $U$ is open and $A$ is closed in $\R^m$. By the proposition above and the definitions in \cite{fed} the current $S \defl \psi_{\#} T$ is an element of $\Flat_n^{\loc}(U)$ resp.\ $\IFlat_n^{\loc}(U)$. Clearly, $\spt(S) = \psi(\spt(T))$ and the result follows with $\varphi \defl \psi^{-1}$.
\end{proof}

In general, it is not possible to take $\alpha = 1$ in the corollary above. To see this, consider the geodesic metric space $G$ which is the Gromov-Hausdorff limit of the so called Laakso graphs, see e.g.\ \cite[p.290]{lang2}. $G$ is doubling and there is no bi-Lipschitz embedding into a Hilbert space as shown in \cite[Theorem 2.3]{lang2}. It is possible to construct a current in $\Norm_1(G)$ with support $G$. Similarly there is a compact geodesic metric space $X$ which is doubling, homeomorphic to $[0,1]^2$, contains an isometric copy of $G$ and is the support of a current in $\Int_2(X)$. The following figure indicates how a homeomorphic image of $X$ in $\R^2$ may look like:

%In general it is not possible to take $\alpha = 1$ in the corollary above. To see this we construct a geodesic metric space $X$ which is doubling, homeomorphic to $[0,1]^2$ and supports an integral current $T \in \Int_2(X)$ yet there is no bi-Lipschitz embedding of $X$ into a Hilbert space. We start with the geodesic metric space $\Gamma$ defined in \cite{lang2} just before Theorem 2.3. Two rectangles are glued to $\Gamma$ and to some of the circles in $\Gamma$ we glue half-spheres with corresponding equators such that we obtain a geodesic metric space $X$ which is homeomorphic to $[0,1]^2$ as indicated in the following figure:

\begin{figure}[h]
\centering
\includegraphics[width=0.75\textwidth, trim = 0cm 0.6cm 0cm 0.45cm]{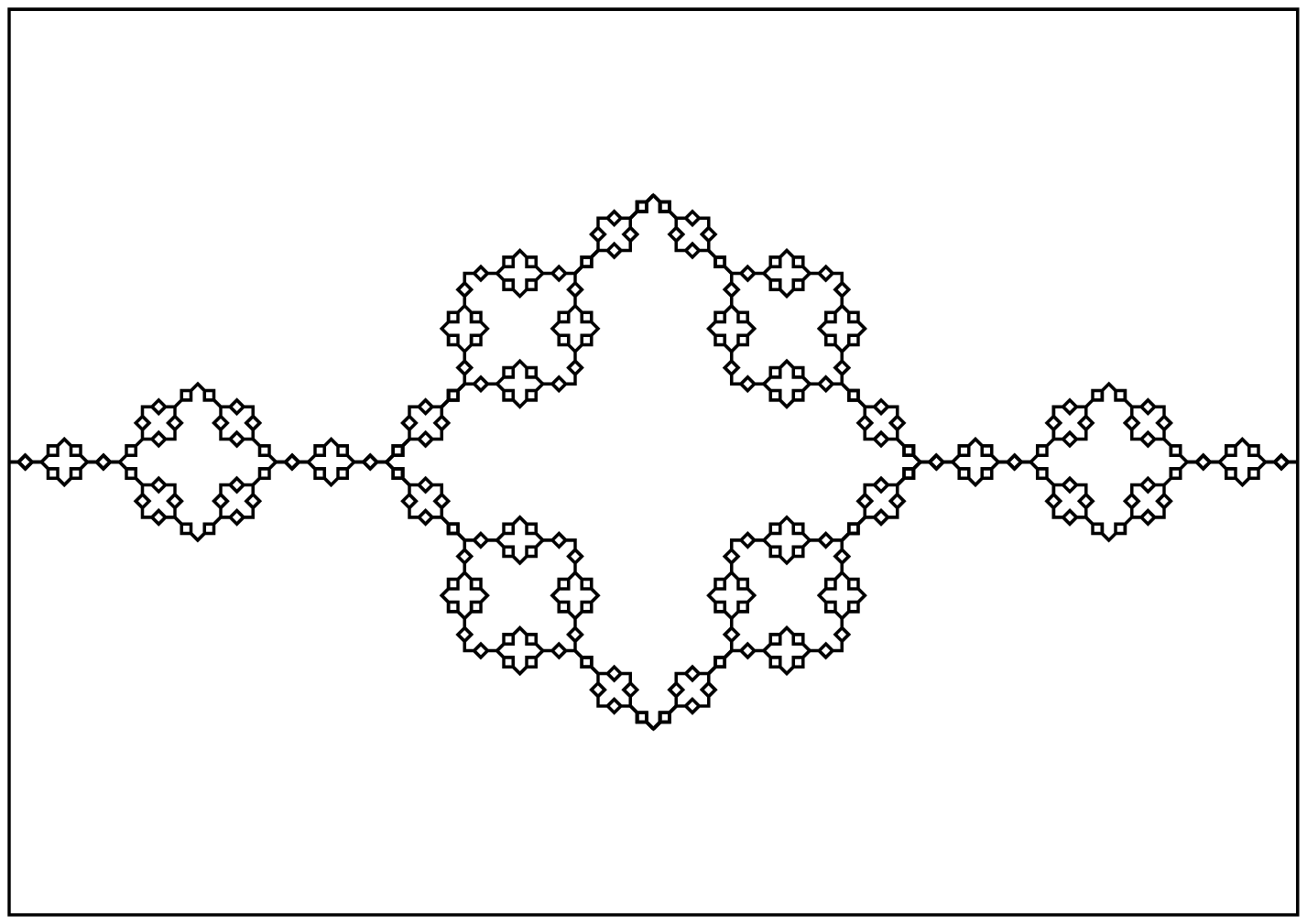}
%\caption{Transverse momentum distributions}
\label{fig1}
\end{figure}

%{\parindent0mm
%The rectangles and half-spheres are oriented in order for their boundaries to cancel along $\Gamma$. This induces an integral current $T \in \Int_2(X)$ with $\spt(T) = X$. One can verify that $\|T\|$ is a doubling measure. Therefore $X$ is doubling and because $X$ contains $\Gamma$ it is not bi-Lipschitz equivalent to a subset of any Hilbert space by \cite[Theorem 2.3]{lang2}.
%}

\subsection{Nonexistence of currents in certain snowflake spaces}

As seen above, a locally normal current in $\Norm_{n,\loc}(X)$ can be extended for certain values, $\alpha > \tfrac{n}{n+1}$, to get a current in $\D_n(X^\alpha)$. The next theorem demonstrates that the bound $\tfrac{n}{n+1}$ is best possible and similar extensions for $\alpha \leq \tfrac{n}{n+1}$ are impossible for all nontrivial currents in $\D_n(X)$.

\begin{thm}
\label{extension2}
Let $X$ be a locally compact metric space and let $n \geq 1$ and $\alpha \leq \tfrac{n}{n+1}$ be constants. Then $\D_n(X^\alpha) = \{0\}$. More generally, for $T \in \D_n(X) \setminus \{0\}$ there is no extension as described in Theorem~\ref{normal} if $\alpha + \beta_1 + \dots + \beta_n \leq n$.
\end{thm}

\begin{proof}
We will give a proof of the first statement because up to some notational changes the second is proved alike. We assume that there is a current $T \in \D_n(X^\alpha)\setminus\{0\}$ and derive a contradiction. Successively, we show that:
\begin{enumerate}
	\item \label{pro1} There is a current $T_1 \in \D_n((\R^n)^\alpha)\setminus\{0\}$ with compact support.
	\item \label{pro2} There is a current $T_2 \in \D_n((\R^n)^\alpha)\cap \Mass_n(\R^n) \setminus \{0\}$ with compact support.
	\item \label{pro3} There is a translation invariant current $T_3 \in \D_n((\R^n)^\alpha)\cap \Mass_{n,\loc}(\R^n) \setminus \{0\}$. Translation invariant means that ${\tau_s}_{\#} T_3 = T_3$ for all $s \in \R^n$, where $\tau_s(x) = x + s$ is the translation by $s$.
	\item \label{pro4} $T_3 = c[\R^n]$ for some $c \neq 0$.
\end{enumerate}
The last point contradicts the fact that $[\R^n] \notin \D_n((\R^n)^\alpha)$ by the counterexample in subsection~\ref{assump}.

\eqref{pro1} Because $T$ restricts to a non-zero current in $\D_n(X)$ there is a $(f,\pi) \in \D^n(X)$ with $T(f,\pi) \neq 0$. The current $T_1 \defl \pi_{\#} (T \lfloor f) \in \D_n(\R^n)$ can be extended to a current in $\D_n((\R^n)^\alpha)$ because $\pi|_{\spt(f)}$ is proper and an element of $\Lip(\spt(f)^\alpha,(\R^n)^\alpha)$. For a function $\sigma \in \Lipc(\R^n)$ with $\sigma = 1$ on a neighborhood of $\pi(\spt(f))$ we get
\[ T_1(\sigma,id) = T_1(f \sigma \circ \pi, id \circ \pi) = T(f,\pi) \neq 0. \]
Hence $T_1 \neq 0$.

\eqref{pro2} Let $\epsilon > 0$. $S_\epsilon \in \D_n((\R^n)^\alpha)$ is defined by
\[ S_\epsilon \defl \frac{1}{\epsilon^n} \int_{[0,\epsilon]^n} {\tau_s}_{\#}T_1 \, d\mathcal{L}^n(s). \]
The map $s \mapsto {\tau_s}_{\#}T_1(f,\pi)$ is continuous for a fixed tuple $(f,\pi) \in \D^n((\R^n)^\alpha)$ by the continuity property of $T_1$. Hence $S_\epsilon$ is a function on $\D^n((\R^n)^\alpha)$ which is multilinear and satisfies the locality condition by definition. To see that $S_\epsilon$ is continuous, Lebesgue's dominated convergence theorem can be used. This is possible because for fixed $C,L \geq 0$ the Arzel\`a-Ascoli theorem and the continuity of $T_1$ imply that the supremum
\begin{equation}
\label{supdef}
%\sup \{ |T_1(f,\pi)| \, : \, (f,\pi) \in \D^n((\R^n)^\alpha), \, \|f\|_\infty \leq C, \, \max\{\Hol^\alpha(f), \max_{1 \leq i \leq n} \Hol^\alpha(\pi^i)\} \leq L \}
\sup \{ |T_1(f,\pi)| \, : \, (f,\pi) \in \Hol_{c,L}^\alpha(\R^n)^{n+1}, \, \|f\|_\infty \leq C \},
\end{equation}
where $\Hol_{c,L}^\alpha(\R^n) \defl \{g \in \Holc^\alpha(\R^n) \, : \, \Hol^\alpha(g) \leq L \}$, is attained and is thus finite. Note that because $\spt(T_1)$ is compact we can assume that each $\pi^i$ in \eqref{supdef} satisfies $\pi^i(x_0) = 0$ for some fixed $x_0 \in \spt(T_1)$ and the support of all functions $f,\pi^1,\dots,\pi^n$ is contained in some compact set depending on $C,L$ and $\spt(T_1)$. So $S_\epsilon$ is indeed a current and
\[ \lim_{\epsilon \rightarrow 0} S_\epsilon(f,\pi) \rightarrow T_1(f,\pi) \]
for every $(f,\pi) \in \D^n((\R^n)^\alpha)$. Hence $S_\epsilon \neq 0$ for $\epsilon$ small enough and we set $T_2 \defl S_\epsilon$ for such an $\epsilon$. Clearly, $T_2$ has compact support. To check that the mass of $T_2$, seen as a current in $\D_n(\R^n)$, is finite, note that
\[ \Mass(T_2) = \sup \{T_2(f,id) \, : \, f\in C_c^\infty(\R^n), \, \|f\|_\infty \leq 1 \}. \]
This follows from the chain rule for currents, \cite[Theorem 2.5]{lang}, and the fact that $C_c^\infty(\R^n)$ is dense in $\D(\R^n)$. Thus, it is enough to consider $T_2(f,id)$ for $f\in \Lipc(\R^n)$ to calculate the mass:
\begin{align*}
T_2(f,id) & = \frac{1}{\epsilon^n}\int_{[0,\epsilon]^n} {\tau_s}_{\#}T_1(f, id) \, d\mathcal{L}^n(s) \\
 & = \frac{1}{\epsilon^n}\int_{[0,\epsilon]^n} T_1(f \circ \tau_s, id + s) \, d\mathcal{L}^n(s) \\
 & = \frac{1}{\epsilon^n}\int_{[0,\epsilon]^n} T_1(f \circ \tau_s, id) \, d\mathcal{L}^n(s) \\
 & = T_1(f_\epsilon, id),
\end{align*}
where $f_\epsilon(x) \defl \frac{1}{\epsilon^n}\int_{[0,\epsilon]^n} f(x + s) \, d\mathcal{L}^n(s)$. In the third line locality and multilinearity of $T_1$ is used. The last equality can be seen by approximating the integral by Riemannian sums and using linearity in the first argument of $T_1$. Choose a sequence $(f_i)_{i\in\N} \subset \Lipc(\R^n)$ with $\|f_i\|_\infty \leq 1$ and $\lim_{i\rightarrow \infty} T_2(f_i,id) = \Mass(T_2)$. Because $\spt(T_2)$ is compact we can assume that for all $i$ the support $\spt(f_i)$ is contained in some fixed compact set $K \subset \R^n$. Clearly $\|(f_i)_\epsilon\|_\infty \leq 1$ for all $i$ and if we can show that $\Lip((f_i)_\epsilon)$ is bounded in $i$, then the Arzel\`a-Ascoli theorem and the continuity of $T_1$ imply that
\[ \Mass(T_2) = T_1(g,id) < \infty \]
for some $g \in \Lipc(\R^n)$. If $\psi$ is any measurable function on $\R^n$ with $|\psi| \leq 1$ almost everywhere, then
\begin{align*}
|\psi_\epsilon(x) - \psi_\epsilon(y)| & = \frac{1}{\epsilon^n} \left| \int_{[0,\epsilon]^n} \psi(x + s) - \psi(y + s) \, d\mathcal{L}^n(s) \right| \\
 & \leq \frac{1}{\epsilon^n} \int_{\tau_x([0,\epsilon]^n) \Delta \tau_y([0,\epsilon]^n)} |\psi(s)| \, d\mathcal{L}^n(s) \\
 & \leq \frac{1}{\epsilon^n}\mathcal{L}^n (\tau_x([0,\epsilon]^n) \Delta \tau_y([0,\epsilon]^n)),
\end{align*}
where $A \Delta B = (A\setminus B) \cup (B\setminus A)$ is the symmetric difference of two sets $A$ and $B$. It is a straight forward calculation to verify that 
\[ \frac{1}{\epsilon^n}\mathcal{L}^n (\tau_x([0,\epsilon]^n) \Delta \tau_y([0,\epsilon]^n)) \leq \frac{L_n}{\epsilon} \|x - y\| \]
for some constant $L_n$ depending only on $n$. Therefore $\Lip(\psi_\epsilon) \leq \tfrac{L_n}{\epsilon}$ and the same holds for the functions $f_i$, which is the remaining part of \eqref{pro2}.

\eqref{pro3} For a current $Z \in \D_n((\R^n)^\alpha)$ with compact support and finite mass we define
\[ \bar Z \defl \int_{\R^n} {\tau_s}_{\#}Z \, d\mathcal{L}^n(s). \]
Because $Z$ has compact support this defines a current in $\D_n((\R^n)^\alpha)$ by the same reasoning as for $T_2$ above. The current $\bar Z$ is apparently translation invariant and has locally finite mass because $Z$ has finite mass and compact support. The crucial part is to find such a $Z$ with $\bar Z \neq 0$. If $f\in \Lipc(\R^n)$ and $s\in \R^n$,
\begin{align}
\nonumber
|{\tau_s}_{\#}Z(f,id) - Z(f,id)| & = |Z(f\circ \tau_s - f,id)| \\
\nonumber
 & \leq \|f\circ \tau_s - f\|_\infty \Mass(Z) \\
\label{lipest}
 & \leq \Lip(f) \|s\| \Mass(Z).
\end{align}
Let $Z$ and $f$ be such that $\spt(Z) \subset B_\epsilon(x_0)$, $\spt(f) \subset B_{\epsilon + \delta}(x_0)$ and 
\[ Z(f,id) > \Lip(f)(2\epsilon + \delta) \frac{n}{n+1}\Mass(Z). \]
Here $B_r(x)$ denotes the closed ball in $\R^n$ with radius $r$ and center $x$. Then \eqref{lipest} implies
\begin{align*}
\bar Z(f,id) %& = \int {\tau_s}_{\#}Z(f,id) \, d\mathcal{L}^n(s) \\
 & = \int_{B_{2\epsilon + \delta}(0)} {\tau_s}_{\#}Z(f,id) \, d\mathcal{L}^n(s) \\
 & \geq \int_{B_{2\epsilon + \delta}(0)} Z(f,id) - \Lip(f) \|s\| \Mass(Z) \, d\mathcal{L}^n(s) \\
 & = \Haus^{n-1}(S^{n-1}) \left( Z(f,id)\frac{(2\epsilon + \delta)^{n}}{n} - \Lip(f)\frac{(2\epsilon + \delta)^{n+1}}{n+1} \Mass(Z) \right) \\
 & > 0,
\end{align*}
which shows that $\bar Z \neq 0$. If $\|f\|_\infty \leq 1$, by altering $f$ outside the ball $B_\epsilon(x_0)$ if necessary, $\delta$ can be assumed to be equal $\Lip(f)^{-1}$. If there is a $Z \in \D_n((\R^n)^\alpha)$ with finite mass and $\spt(Z) \subset B_\epsilon(x_0)$ and an $f \in \Lipc(\R^n)$ with $\|f\|_\infty \leq 1$ such that
\begin{equation}
\label{import}
Z(f,id) > (2\epsilon\Lip(f) + 1)\frac{n}{n+1}\Mass(Z),
\end{equation}
the $T_3$ we look for can be constructed. To find such a $Z$ fix $\rho \in (\tfrac{n}{n+1},1)$ and let $f$ be an element of $\Lipc(\R^n)$ with $\|f\|_\infty \leq 1$ and $T_2(f,id) > \rho \Mass(T_2)$. Choose $\epsilon > 0$ such that $\rho \geq (2\epsilon\Lip(f) + 1) \frac{n}{n+1}$ and a Lipschitz partition of unity $f_1,\dots,f_N$ in $\R^n$ for $\spt(T_2)$ such that $\spt(f_i) \subset B_{\epsilon}(x_i)$ for some $x_i \in \R^n$, $i=1,\dots,N$. Now,
\[ \sum_{i=1}^N (T_2 \lfloor f_i) (f,id) = T_2(f,id) > \rho \Mass(T_2) = \sum_{i=1}^N \rho \Mass(T_2 \lfloor f_i). \]
Hence there is at least one $i$ such that
\[ (T_2 \lfloor f_i) (f,id) > \rho \Mass(T_2 \lfloor f_i) \geq (2\epsilon\Lip(f) + 1)\frac{n}{n+1} \Mass(T_2 \lfloor f_i). \]
Clearly $T_2 \lfloor f_i$ is a current in $\D_n((\R^n)^\alpha)$ with $\spt(T_2 \lfloor f_i) \subset B_{\epsilon}(x_i)$ and finite mass for which \eqref{import} holds.

\eqref{pro4} By construction $\|T_3\|$ is a nontrivial, translation invariant Radon measure on $\R^n$. See \cite[Chapter 4]{lang} for the definition and properties of this measure. Thus $\|T_3\| = C\mathcal{L}^n$ for a constant $C > 0$. We set $T_3(\chi_{[0,1]^n},id) \defr c$. By the linearity of $T_3$ in the first argument, the inequality $|T_3(f,id)| \leq C\int|f| \, d\mathcal {L}^n$ and Lebesgue's dominated convergence theorem we conclude that
\[ T_3(f,id) = c \int_{\R^n} f(x) \, d\mathcal {L}^n(x) \]
for all $f\in \mathcal B_c^\infty(\R^n)$. This also implies that $C = |c|$. If $\pi \in C^{1,1}(\R^n,\R^n)$, the identity
$ T_3(f,\pi) = c[\R^n](f,\pi) $
holds by the chain rule, \cite[Theorem 2.5]{lang}. And finally $T_3 = c[\R^n]$ by approximating the Lipschitz functions with smooth ones.
\end{proof}

An immediate consequence of this result is that $\D_n(X) = \{0\}$ for some $n \geq 1$ if $X$ is bi-Lipschitz equivalent to a locally compact metric space $Y$ for which
\[ d(x,z)^{\frac{n+1}{n}} \leq d(x,y)^{\frac{n+1}{n}} + d(y,z)^{\frac{n+1}{n}} \]
holds for all $x,y,z \in Y$. For instance, this is true for all $n \geq 1$ if $Y$ is an ultrametric space, i.e.
\[ d(x,z) \leq \max\{d(x,y), d(y,z)\} \]
holds for all $x,y,z \in Y$.
% then $\D_n(X) = \{0\}$ for all $n \geq 1$ and any $X$ bi-Lipschitz equivalent to $Y$.

\bibliographystyle{bibstyle}
\bibliography{refs}

\end{document}